\documentclass{article}

\usepackage{arxiv}

\usepackage[utf8]{inputenc}
\usepackage[T1]{fontenc}
\usepackage{amsmath} 
\usepackage{amsthm}
\usepackage{amssymb}
\usepackage{amsfonts}  
\usepackage[normalem]{ulem}
\usepackage{algorithm}
\usepackage{algpseudocodex}
\usepackage{hyperref}       
\usepackage{url}            
\usepackage{booktabs}    
\usepackage{placeins}
\usepackage{nicefrac}      
\usepackage{microtype}      
\usepackage{cleveref}       
\usepackage{lipsum}        
\usepackage{graphicx}
\usepackage[numbers]{natbib}
\usepackage{doi}
\usepackage{todonotes}
\usepackage{caption}
\usepackage{subcaption}
\usepackage{siunitx}

\usepackage{tikz}
\usetikzlibrary{calc}
\usetikzlibrary{graphs, graphs.standard}
\usetikzlibrary{decorations.pathreplacing, calligraphy, matrix}

\theoremstyle{definition}
\newtheorem{proposition}{Proposition}

\newtheorem{observation}{Observation}

\DeclareMathOperator{\Ord}{Ord}

\newcommand{\gmin}{g_{\min}}

\usepackage{titlesec}
\titleformat{\paragraph}
{\normalfont\normalsize\bfseries}{\theparagraph}{1em}{}
\titlespacing*{\paragraph}
{0pt}{1.25ex plus 0.5ex minus .2ex}{0.5ex plus .2ex}

\title{New small regular graphs of given girth: the cage problem and beyond}

\author{ %
\href{https://orcid.org/0000-0003-2023-8887}{\includegraphics[scale=0.06]{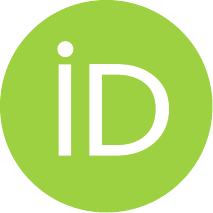}\hspace{1mm}Geoffrey Exoo} \\
	Department of Mathematics and Computer Science \\
	 Indiana State University\\
	 Terre Haute, IN 47809, USA \\
	\texttt{ge@cs.indstate.edu} \\
	\And
	\href{https://orcid.org/0000-0001-8984-2463}{\includegraphics[scale=0.06]{orcid.pdf}\hspace{1mm}Jan Goedgebeur} \\
	Department of Computer Science\\
	KU Leuven Kulak\\
    8500 Kortrijk, Belgium \\
	\texttt{jan.goedgebeur@kuleuven.be} \\
    \And
	\href{https://orcid.org/0000-0002-5256-1921}{\includegraphics[scale=0.06]{orcid.pdf}\hspace{1mm}Jorik Jooken} \\
	Department of Computer Science\\
	KU Leuven Kulak\\
    8500 Kortrijk, Belgium \\
	\texttt{jorik.jooken@kuleuven.be} \\
    \And
	\href{https://orcid.org/0009-0004-9873-8966}{\includegraphics[scale=0.06]{orcid.pdf}\hspace{1mm}Louis Stubbe} \\
	Department of Computer Science\\
	KU Leuven campus Ghent\\
    9000 Ghent, Belgium \\
	\texttt{louis.stubbe@kuleuven.be} \\
    \And
	\href{https://orcid.org/0009-0003-7686-1982}{\includegraphics[scale=0.06]{orcid.pdf}\hspace{1mm}Tibo Van~den~Eede} \\
	Department of Computer Science\\
	KU Leuven Kulak\\
    8500 Kortrijk, Belgium \\
	\texttt{tibo.vandeneede@kuleuven.be} \\
}

\renewcommand{\shorttitle}{New small regular graphs of given girth}

\hypersetup{
pdftitle={New small regular graphs of given girth: the cage problem and beyond},
pdfsubject={cs.DM},
pdfauthor={Geoffrey Exoo, Jan Goedgebeur, Jorik Jooken, Louis Stubbe, Tibo Van~den~Eede},
pdfkeywords={Cage problem, Extremal graph theory, Lifts of voltage graphs, Graph generation algorithm},
}

\begin{document}
\maketitle

\begin{abstract}
The cage problem concerns finding $(k,g)$-graphs, which are $k$-regular graphs with girth $g$, of the smallest possible number of vertices. The central goal is to determine $n(k,g)$, the minimum order of such a graph, and to identify corresponding extremal graphs. In this paper, we study the cage problem and several of its variants from a computational perspective. Four complementary graph generation algorithms are developed based on exhaustive generation of lifts, a tabu search heuristic, a hill climbing heuristic and excision techniques. Using these methods, we establish new upper bounds for eleven cases of the classical cage problem: $n(3,16) \leq 936$, $n(3,17) \leq 2048$, $n(4,9) \leq 270$, $n(4,10) \leq 320$, $n(4,11) \leq 713$, $n(5,9) \leq 1056$, $n(6,11) \leq 7783$, $n(8,7) \leq 774$, $n(10,7) \leq 1608$, $n(12,7) \leq 2890$ and $n(14,7) \leq 4716$. Notably, our results improve upon several of the best-known bounds, some of which have stood unchanged for 22 years. Moreover, the improvement for $n(4,10)$, from the longstanding upper bound of 384 down to 320, is surprising and constitutes a substantial improvement.

While the main focus is on the cage problem, we also adapted our algorithms for variants of the cage problem that received attention in the literature. For these variants, additional improvements are obtained, further narrowing the gaps between known lower and upper bounds.

Note: The previous arXiv version of this paper established $n(5,9) \leq 1116$. A subsequent search using similar voltage-lift methods yielded a construction further improving the bound to $n(5,9) \leq 1056$.
\end{abstract}


\section{Introduction}\label{sec:introduction}

A \emph{$(k,g)$-graph} is a $k$-regular graph with girth $g$, i.e., a graph in which every vertex has degree $k$ and the length of the shortest cycle is $g$. For all integers $k \geq 2, g \geq 3$ it is known that a $(k,g)$-graph exists~\cite{sachs1963}. A natural question is then to wonder how small such a graph can be, making these graphs central objects of study in extremal graph theory. The cage problem is concerned with exactly this question. More precisely, a \emph{$(k,g)$-cage} is a $(k,g)$-graph with the smallest possible number of vertices, denoted by $n(k,g)$, and the cage problem is concerned with determining $n(k,g)$ and the corresponding $(k,g)$-cages. A graph that is a $(k,g)$-cage for some $k$ and $g$ will simply be referred to as a cage.

Cages are not only of theoretical interest. Their sparse and often highly symmetric structure makes them attractive for practical applications in areas such as communication networks~\citep{DonettiLuca2006Onte, RoigPedroJuan2023EDCO}, error-correcting codes~\citep{malema2007high, YangHongna2024Erga}, and cryptography~\citep{CsirmazLászló2019Ssol}. For instance, in the context of network design, cages provide topologies that support efficient and uniform information propagation, avoid overburdened central nodes and exhibit robustness against failures. 

Despite the seemingly simple description, finding $(k,g)$-cages is notoriously difficult. For most parameter pairs $(k,g)$, the exact value of $n(k,g)$ is unknown. In those cases, research has focused on finding good lower and upper bounds for $n(k,g)$ instead. Several such bounds are known in the literature. Let us introduce one classical lower bound for $n(k,g)$, namely the Moore bound $M(k,g)$. Let $G$ be a $(k,g)$-graph and let $u$ be any vertex in $V(G)$. Then there cannot be an edge between any two vertices $v$ and $w$ if the distance between $u$ and $v$ plus the distance between $u$ and $w$ is strictly smaller than $g-1$. If $g$ is odd, one can consider vertices that are close to a vertex $u$, whereas if $g$ is even, one can consider vertices that are close to an edge $uv$. In both cases, there must be enough vertices to ensure that at least one cycle of length $g$ can be present and this argument gives rise to the Moore bound $M(k,g)$:
\[
n(k,g) \;\ge\; M(k,g)= 
\begin{cases}
1 + k \displaystyle\sum_{i=0}^{(g-3)/2} (k-1)^i, & \text{if $g$ is odd}, \\[10pt]
2 \displaystyle\sum_{i=0}^{(g-2)/2} (k-1)^i, & \text{if $g$ is even}.
\end{cases}
\]

A $(k,g)$-graph with order $M(k,g)$ is called a Moore graph. Moore graphs can only exist for $g \in \{3,4,5,6,8,12\}$ and for many $(k,g)$ pairs, the best lower bound on $n(k,g)$ is given by $M(k,g)+1$ or $M(k,g)+2$ \citep{cagesurvey}. For small $k$ and $g$, however, better bounds have often been established using clever computational approaches \citep{jooken2025}. An overview of lower and upper bounds can be found in the Dynamic Cage Survey \citep{cagesurvey}, while some recent improvements (including the improvements from the current paper) are available on \href{https://houseofgraphs.org/meta-directory/cages}{The House of Graphs} \citep{houseofgraphs} at \url{https://houseofgraphs.org/meta-directory/cages}.

Given the difficulty of improving the bounds on $n(k,g)$, researchers have introduced and explored several variants of the cage problem, both for their independent interest and for the insights they may offer into the classical case. In this paper we will therefore also study some of the main variants. One such variant involves the study of \emph{$(k,g)$-spectra}~\citep{eze2025spectra}: the set of all orders $n$ for which a $(k,g)$-graph exists. While determining the spectrum is at least as hard as solving the cage problem itself, partial knowledge can be helpful in excluding the existence of certain graphs. One such example is the hypothetical $(57,5)$-Moore graph of order 3250 (the existence of which is a famous open problem~\cite{DALFO20191}). More specifically, if the existence of some $(57,5)$-graph of order $n>3250$ can be excluded, but a construction algorithm developed in the context of $(k,g)$-spectra could create such a graph starting from a candidate $(57,5)$-Moore graph, then the nonexistence of such a candidate can be established indirectly.

A second and third variant are vertex-girth-regular and edge-girth-regular graphs, denoted respectively as $(k,g,\lambda_v)$- and $(k,g,\lambda_e)$-graphs~\citep{goedgebeur2024egrgraphs, JAJCAY201870, jajcay2024vgr} or $vgr$- and $egr$-graphs. A graph is \emph{vertex-girth-regular} if it is a $(k,g)$-graph in which every vertex lies on exactly $\lambda_v$ cycles of length $g$ for some integer $\lambda_v$. Similarly, a graph is \emph{edge-girth-regular} if it is a $(k,g)$-graph in which every edge lies on exactly $\lambda_e$ cycles of length $g$ for some integer $\lambda_e$. The goal is then to find such graphs (if they exist) of the smallest order, denoted as $n(k,g,\lambda_v)$ and $n(k,g,\lambda_e)$ for $vgr$- and $egr$-graphs, respectively. This structural constraint is inspired by the observation that many known cages exhibit a high degree of regularity in their cycle distribution and therefore many $(k,g)$-cages are either vertex-girth-regular or edge-girth-regular graphs. As such, only searching for $vgr$- and $egr$-graphs can act as a heuristic, narrowing the search space.

A fourth and final variant studied in this paper involves $(k,g)$-graphs that contain no cycles of length $g+1$, denoted as $(k,g,\uline{g+1})$-graphs~\citep{eze2024kggraphsg1cycles}. This again leads to an extremal problem where the goal is to find the minimum order of such a graph, denoted by $n(k,g,\uline{g+1})$. These graphs are particularly relevant when $g$ is odd, due to the inequality $n(k,2s+4) \leq 2n(k,2s+1,\uline{2s+2})$~\cite{eze2024kggraphsg1cycles} that connects these graphs with the cage problem.

In this paper, we establish several new upper bounds for the classical cage problem. Specifically, we improve the best known values to  
\[
n(3,16) \leq 936,\quad n(3,17) \leq 2048,\quad n(4,9) \leq 270,\quad n(4,10) \leq 320, \quad n(4,11) \leq 713, \quad n(5,9) \leq 1056,
\]  
\[
n(6,11) \leq 7783,\quad 
n(8,7) \leq 774, \quad 
n(10,7) \leq 1608,\quad 
n(12,7) \leq 2890, \quad 
n(14,7) \leq 4716.
\]

Most of these results represent significant progress compared to the previous bounds, in particular:
\[
n(3,16) \leq 960 \cite{exoovoltage2004},\quad n(3,17) \leq 2176\cite{exoo2019new},\quad n(4,9) \leq 275\cite{cagesurvey},\quad n(4,10) \leq 384\cite{cagesurvey},\quad n(5,9) \leq 1152\cite{Exoo2025}.
\]  
We remark that some of these bounds had stood unchanged for 22 years. Moreover, the reduction of the upper bound for $n(4,10)$ from 384 to 320 marks a substantial improvement. To the best of our knowledge, the previous best upper bounds on $n(4,11)$ and $n(6,11)$ were the bounds $n(4,12)-1$ and $n(6,12)-1$, respectively. These bounds are obtained from the inequality $n(k, g) < n(k, g + 1)$, which holds for every $k\ge2$ and $g\ge3$ \cite{sauer1967extremal}. Thus, the improvement on $n(4,11)$ and $n(6,11)$ marks the first non-trivial upper bound for these cases. 

In addition to these classical cases, we also obtained multiple improvements for several variants of the cage problem. Specifically, we found 21, 29, and 6 new upper bounds for $egr$-, $vgr$-, and $(k,g,\uline{g+1})$-graphs, respectively. Among these, 2 bounds for $egr$-graphs and 7 bounds for $vgr$-graphs were tight, i.e., equal to the corresponding lower bounds. Furthermore, we resolved 34 previously open cases for $(k,g)$-spectra by constructing graphs for these unresolved orders.

The rest of this paper is structured as follows. In Section~\ref{sec:prelim} we introduce preliminaries regarding groups, graphs and lifts of voltage graphs. Next, in Section~\ref{sec:algorithms} we present the different algorithms that were developed and implemented to construct the $(k,g)$-graphs that allowed us to improve several of the previously best bounds. We present a detailed overview of these results in Section~\ref{sec:results}. Finally, we close this paper in Section~\ref{sec:conclusion} with some concluding remarks.

\section{Preliminaries}\label{sec:prelim}
Throughout this paper, all groups are assumed to be finite. For a group $\Gamma$, we denote its identity element by $0_\Gamma$ and the inverse of an element $r \in \Gamma$ by $r^{-1}$. $C_n$ and $D_n$ denote the cyclic group with $n$ and the dihedral group with $2n$ elements, respectively. For an element $r \in \Gamma$, we say the \emph{order} of $r$, denoted by $\Ord(r)$, is the smallest integer $\ell$ such that $r^\ell=0_\Gamma$, where $r^\ell$ denotes the product of $\ell$ copies of $r$.

A \emph{digraph} $D$ is a pair $(V, A)$ where $V$ is a finite set of vertices and $A$ is a multiset 
of ordered pairs of vertices, called \emph{arcs}. Thus multiple arcs between the same ordered pair of 
vertices are allowed, and loops (arcs of the form $(u,u)$) are also permitted.

A \emph{walk} $W$ from vertex $v_0$ to $v_n$ in digraph $D = (V, A)$ consists of an alternating sequence of vertices $v_i \in V$ and arcs $a_i \in A$:
\[
W = v_0, a_0, v_1, a_1, \dots , a_{n-1}, v_n, 
\]
with $a_i$ an arc between $v_i$ and $v_{i+1}$. A walk is \emph{open} if $v_0 \neq v_n$ and 
\emph{closed} otherwise.

A \emph{graph} $G$ is a pair $(V,E)$ where $V$ is a finite set of vertices and $E$ is a multiset of 
\emph{edges}. Each edge is either:
\begin{itemize}
    \item a \emph{full edge} $\{u,v\}$, a multiset containing two (not necessarily distinct) vertices 
    $u,v \in V$, or
    \item a \emph{semi-edge} $\{u\}$, consisting of a single vertex.
\end{itemize}
Hence, graphs may contain loops, parallel edges and semi-edges. 

We now also introduce the notion of \emph{darts} to be able to naturally associate a digraph with each graph. For every full edge $e=\{u,v\}$, we arbitrarily designate one of the two possible orientations as the positive dart $e^+=(u,v)$ and the other as the negative dart $e^-=(v,u)$. These darts are inverses of each other, denoted $(e^+)^{-1}=e^-$ and $(e^-)^{-1}=e^+$. A loop $\{u,u\}$ gives rise to two darts $(u,u)$, each inverse of the other. A semi-edge $\{u\}$ produces a single dart $(u,u)$, which is its own inverse.

Let $D=(V,A)$ be a digraph and let $\Gamma$ be a group. A \emph{voltage assignment} on $D$ is a function $\alpha : A \to \Gamma$. The \emph{lift} (or \emph{derived regular cover}) $D^\alpha$ is the digraph with vertex set
\[
V(D^\alpha) = \{\, u^s \mid u \in V, \ s \in \Gamma \,\},
\]
and arc set $A(D^\alpha)$ defined as follows: for every arc $a=(u,v) \in A$ and every $s \in \Gamma$, there is an arc
\[
(u^s, \, v^{s \cdot \alpha(a)}) \in A(D^\alpha).
\]

To apply the voltage construction to an undirected graph $G=(V,E)$ we first replace $G$ with its dart digraph $D(G)=(V,A)$ as described above. A 
voltage assignment is then a function $\alpha : A \to \Gamma$. In order to be able to interpret the resulting lift $D(G)^\alpha$ as a graph (by pairing two opposing arcs into a single edge), we assume that the following conditions are met:
\begin{itemize}
\item $\alpha(d^{-1}) = \alpha(d)^{-1}$ for every dart $d \in A$. Note that $d^{-1}$ denotes the structural inverse of dart $d$, while $\alpha(d)^{-1}$ denotes the algebraic inverse of group element $\alpha(d)$.
\item $\alpha(d) \neq 0_{\Gamma}$ for every dart $d \in A$ that originates from a semi-edge.
\end{itemize}

\begin{figure}
\begin{subfigure}[t]{0.25\textwidth}
	\centering
	\begin{tikzpicture}[
		node/.style={circle, draw, minimum size=0.55cm, inner sep=0pt},
		label/.style={midway, fill=white, font=\sffamily\small}
		]
		
		\node[node, fill=red!50] (Ao) at (-3.50,-0.75) {r};
		\node[node, fill=blue!50] (Bo) at (-2.50,-0.75) {b};
		
		\draw[thick] (Ao) to (Bo);
		
		\draw[thick] (Ao) to[out=135, in=225, looseness=6] (Ao);
		
		\draw[thick] (Bo) to[out=45, in=315, looseness=6] (Bo);	
	\end{tikzpicture}
	\caption{The dumbbell graph}
\end{subfigure}
\begin{subfigure}[t]{0.4\textwidth}
	\centering
	\begin{tikzpicture}[
		node/.style={circle, draw, minimum size=0.55cm, inner sep=0pt},
		label/.style={midway, fill=white, font=\sffamily\small}
		]

        \useasboundingbox (0,-2) rectangle (2.5,0.5);
		
		\node[node, fill=red!50] (A) at (0.50,-0.75) {r};
		\node[node, fill=blue!50] (B) at (2,-0.75) {b};
		
		\draw[->, thick, >=stealth] (A) to[out=30, in=150] node[label] {0} (B);
		\draw[->, thick, >=stealth] (B) to[out=-150, in=-30] node[label] {0} (A);
		
		\draw[->, thick, >=stealth] (A) to[out=135, in=225, looseness=7] node[label] {2} (A);
		\draw[->, thick, >=stealth] (A) to[out=245, in=115, looseness=17] node[label] {3} (A);
		
		\draw[->, thick, >=stealth] (B) to[out=45, in=315, looseness=7] node[label] {1} (B);
		\draw[->, thick, >=stealth] (B) to[out=295, in=65, looseness=17] node[label] {4} (B);	
	\end{tikzpicture}
	\caption{The dumbbell graph with darts and voltages}
\end{subfigure}
\begin{subfigure}[t]{0.33\textwidth}
	\centering
	\begin{tikzpicture}[
		node/.style={circle, draw, minimum size=0.55cm, inner sep=0pt},
		label/.style={midway, fill=white, font=\sffamily\small}
		]
		
		\node[node, fill=blue!50] (B1) at (4.0,0.5) {$b^4$};
		\node[node, fill=red!50] (R1) at (5.25,0.25) {$r^4$};
		
		\node[node, fill=blue!50] (B2) at (6.25,2.25) {$b^0$};
		\node[node, fill=red!50] (R2) at (6.25,1) {$r^0$};
		
		\node[node, fill=blue!50] (B3) at (8.5,0.5) {$b^1$};
		\node[node, fill=red!50] (R3) at (7.25,0.25) {$r^1$};
		
		\node[node, fill=blue!50] (B4) at (4.75,-1.75) {$b^3$};
		\node[node, fill=red!50] (R4) at (5.75,-0.75) {$r^3$};
		
		\node[node, fill=blue!50] (B5) at (7.75,-1.75) {$b^2$};
		\node[node, fill=red!50] (R5) at (6.75,-0.75) {$r^2$};
		
		\foreach \x/\y in {B1/R1,B2/R2,B3/R3,B4/R4,B5/R5}{
			\draw[thick] (\x) -- (\y);
		}
		
		\foreach \x/\y in {B1/B2,B2/B3,B3/B5,B5/B4,B4/B1}{
			\draw[thick] (\x) -- (\y);
		}
		
		\foreach \x/\y in {R1/R3,R1/R5,R2/R4,R2/R5,R3/R4}{
			\draw[thick] (\x) -- (\y);
		}	
	\end{tikzpicture}
	\caption{The resulting lift}
\end{subfigure}
	\caption{\label{fig:liftDumbbell} A lift of the dumbbell graph with the group $C_5$ resulting in the Petersen graph. Opposing arcs in the lift were replaced with a single edge.}
\end{figure}
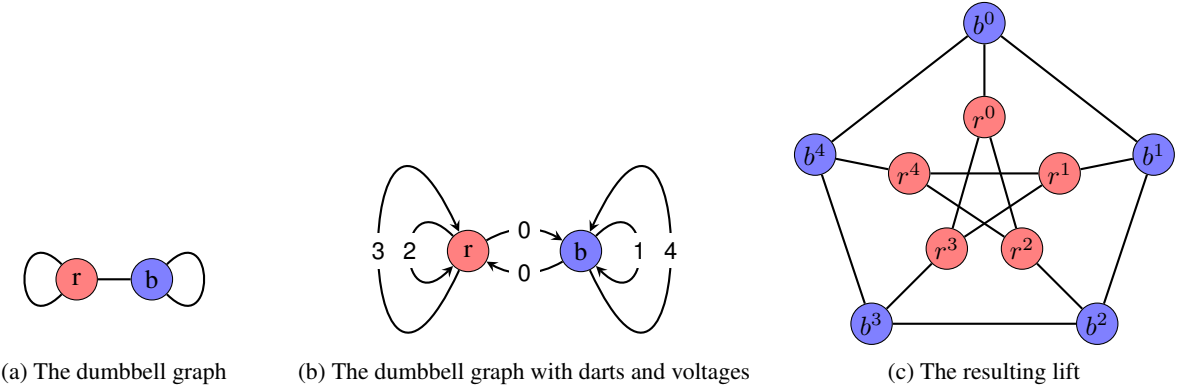

We will refer to $G$ as the \emph{base graph} corresponding to the lift $D(G)^\alpha$. We will often also use $G^\alpha$ to refer to the lift $D(G)^\alpha$. In the rest of this paper, we assume that voltage assignments are chosen so that lifts can be regarded as graphs. We adopt the convention that the girth of a graph containing a loop or semi-edge is $1$, and the girth of a graph containing parallel edges is at most $2$. Since we are interested in $(k,g)$-graphs with $g \geq 3$, this ensures that we only need to consider lifts that produce graphs without loops, parallel edges, or semi-edges. (However, the base graph itself may still contain such features.) Fig.~\ref{fig:liftDumbbell} illustrates an example of lifting the dumbbell graph into the Petersen graph.

\subsection{Properties of lifts} \label{subsubsec:properties}
Lifted graphs inherit several structural properties from their base graph, their group and their chosen voltage assignment. These properties are exploited by the algorithms discussed in the next section. Throughout this subsection, we let $G$ be a graph, $\Gamma$ a group, and $\alpha$ a voltage assignment with corresponding lift $G^\alpha$. The first two observations follow directly from the definition of lifts: 
\begin{observation} \label{ob:regular}
	$G$ is $k$-regular if and only if $G^\alpha$ is $k$-regular.
\end{observation}

\begin{observation} \label{ob:connected}
	If $G$ is not connected, then $G^\alpha$ is not connected.
\end{observation}
Connectedness of the base graph is necessary, but not sufficient, for the lift to be connected. 

It is possible to calculate the girth of a lift without explicitly constructing this lift by finding the shortest, closed, non-reversing walk with net voltage $0_\Gamma$ in the digraph associated with the base graph. The \emph{net voltage} of a walk $W = v_0, a_0, v_1, a_1, \dots , a_{n-1}, v_n$ in the digraph associated with the base graph $G$ with voltage assignment $\alpha$ is the product of all voltages on this walk: $\alpha(a_0) \cdot  \alpha(a_1) \cdot \ldots \cdot \alpha(a_{n-1})$. A walk is \emph{non-reversing} if after using an arc $(u,v)$, it is not immediately followed by using the arc $(v,u)$.

\begin{proposition}[{\cite[Lemma 2.1]{exoo2011voltagegraphlifts}}]\label{prop:walkgirth}
	The girth of $G^\alpha$ is equal to the length of the shortest, closed, non-reversing walk with net voltage $0_\Gamma$ in the digraph associated with $G$.
\end{proposition}

A final set of properties revolves around base graphs, groups and voltage assignments that yield isomorphic lifts. In this context, the following two properties are fairly obvious:

\begin{observation} \label{ob:groupiso}
	For any automorphism $\phi: \Gamma \to \Gamma$, the lift $G^{\alpha^\phi}$ constructed using the voltage assignment $\alpha^\phi$, defined by
	\[
	\forall d \in A(D(G)) : \alpha^\phi(d) = \phi(\alpha(d)),
	\]
	is isomorphic to the original lift $G^\alpha$.
\end{observation}

\begin{observation} \label{ob:graphiso}
Suppose $\phi: E(G) \to E(G)$ is an edge automorphism of $G$, i.e., a bijection of edges such that if two edges in $G$ share a vertex, their images under $\phi$ also share a vertex. Then $\phi$ induces a natural bijection of darts, denoted $\psi: A(D(G)) \to A(D(G))$, by mapping each dart of an edge $e$ to the corresponding dart of the edge $\phi(e)$.  

Define a new voltage assignment $\alpha^\phi : A(D(G)) \to \Gamma$ by
\[
\alpha^\phi(d) := \alpha(\psi(d)) \quad \text{for all darts } d \in A(D(G)).
\]

Then the lift $G^{\alpha^\phi}$ is isomorphic to $G^\alpha$.
\end{observation}

A less obvious source of isomorphic lifts are spanning trees in the base graph. The following proposition greatly reduces the number of options to explore when one wants to compute all graphs that can arise as the lift of a given base graph.
\begin{proposition}[{\cite[p. 91]{topologicalGraphTheory}}]\label{prop:tree}
	Let $T$ be a spanning tree of $G$. There exists a voltage assignment $\alpha'$ with $\alpha'(d) = 0_\Gamma$ for every dart $d$ corresponding with an edge of $T$ such that $G^{\alpha'}$ is isomorphic with $G^\alpha$.
\end{proposition}

\section{Algorithms for generating regular graphs of given girth}\label{sec:algorithms}

The construction of lifts requires traversing a combinatorial search space of possible group–graph combinations and voltage assignments. To manage this complexity, we employ a layered approach that progressively narrows the search space. First, we generate admissible base graphs and groups as described in Subsection~\ref{subsec:selectionbasegraphsgroups}, after which structure- and assignment-based exclusion rules prune infeasible or redundant voltage assignments as shown in Subsection~\ref{subsec:selectionvoltage}. The efficient filtering procedures of Subsection~\ref{subsec:efficientfiltering} then remove isomorphic or unsuitable lifts. Subsection~\ref{subsec:iteratingvoltages} describes how the remaining search space is explored using two complementary strategies: an exhaustive backtracking algorithm that guarantees completeness, and a heuristic tabu search that trades completeness for scalability. In Subsection~\ref{subsec:hill climbing}, we take a different approach and discuss a hill climbing algorithm that simultaneously searches for a suitable base graph and voltage assignment for a given group.

Finally, instead of constructing a $(k,g)$-graph through a lift, Subsection~\ref{subsec:excision} shows how vertices can be excised from an existing small $(k,g)$-graph and how edges can be added again to obtain a smaller $k$-regular graph of large girth.

\subsection{Selection and generation of groups and base graphs} \label{subsec:selectionbasegraphsgroups}
As Subsection~\ref{subsubsec:properties} imposes no restrictions on the usable groups, any collection of pairwise non-isomorphic groups can be utilized. An exhaustive list with all groups of order $n < 1024$ can be obtained using GAP~\cite{GAP4}.

The selection of base graphs is more interesting, since these do have restrictions imposed on them. More specifically, due to Observations~\ref{ob:regular}~and~\ref{ob:connected}, the base graphs must be regular and connected, respectively (since all graphs of interest to the problems mentioned in the introduction require regularity and connectedness). For a fixed order $n$ and regularity $k$, an exhaustive collection of pairwise non-isomorphic graphs with potentially parallel edges, loops and semi-edges can be generated by adapting \texttt{multigraph}. This is a non-published graph generator developed by Gunnar Brinkmann based on the generation algorithm described in \cite{brinkmann2013generation}.

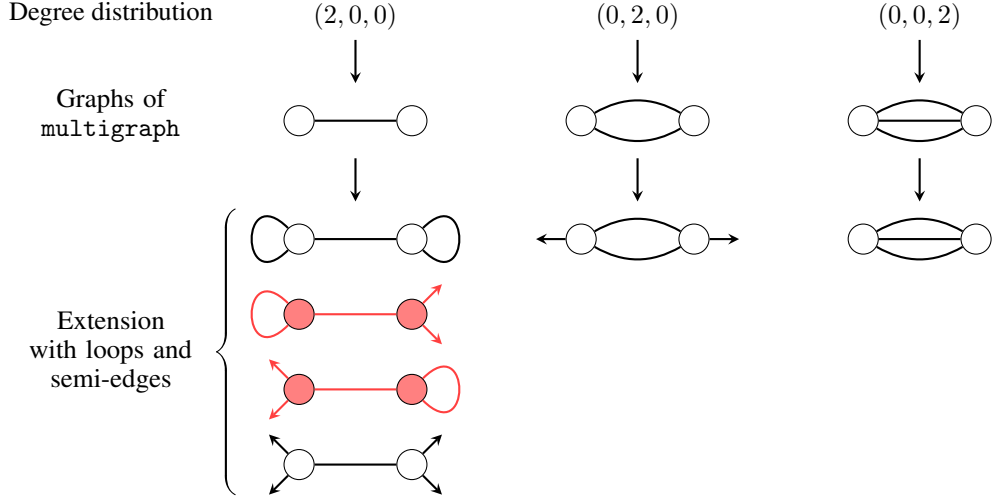
\begin{figure}
	\begin{tikzpicture}[
		node/.style={circle, draw, minimum size=0.4cm},
		label/.style={midway, fill=white, font=\sffamily\small}
		]

		\node[node, draw=none] (1) at (-2.5,0.99) {Degree distribution};
		\node[node, draw=none,align=center,text width=2cm] (1) at (-2.5,-0.35) {Graphs of \texttt{multigraph}};
		\node[node, draw=none,align=center,text width=2.5cm] (1) at (-2.5,-3.5) {Extension with loops and semi-edges};
		\draw [decorate,decoration = {calligraphic brace, mirror, amplitude=7pt},thick] (-0.85,-1.6) --  (-0.85,-5.4);
		
		\draw[<-, thick, >=stealth] (0.75,0.075) -- (0.75,0.65) node[above] {$(2, 0, 0)$} ;
		\draw[<-, thick, >=stealth] (4.5,0.075) -- (4.5,0.65) node[above] {$(0, 2, 0)$};
		\draw[<-, thick, >=stealth] (8.25,0.075) -- (8.25,0.65) node[above] {$(0, 0, 2)$};
		
		\node[node] (1) at (0,-0.425) {};
		\node[node] (2) at (1.5,-0.425) {};
		\draw[thick] (1) to (2);
		
		\draw[->, thick, >=stealth] (0.75,-0.925) -- (0.75,-1.5);
		
		\node[node] (3) at (3.75,-0.425) {};
		\node[node] (4) at (5.25,-0.425) {};
		\draw[thick] (3) to[out=30, in=150] (4);
		\draw[thick] (4) to[out=-150, in=-30] (3);
		
		\draw[->, thick, >=stealth] (4.5,-0.925) -- (4.5,-1.5);
		
		\node[node] (5) at (7.5,-0.425) {};
		\node[node] (6) at (9,-0.425) {};
		\draw[thick] (5) to[out=30, in=150] (6);
		\draw[thick] (5) to (6);
		\draw[thick] (6) to[out=-150, in=-30] (5);
		
		\draw[->, thick, >=stealth] (8.25,-0.925) -- (8.25,-1.5);

		\node[node] (7) at (0,-2) {};
		\node[node] (8) at (1.5,-2) {};
		\draw[thick] (7) to (8);
		\draw[thick] (7) to[out=135, in=225, looseness=8] (7);
		\draw[thick] (8) to[out=45, in=315, looseness=8] (8);
		
		\node[node, fill=red!50] (9) at (0,-3) {};
		\node[node, fill=red!50] (10) at (1.5,-3) {};
		\draw[thick, red!75] (9) to (10);
		\draw[thick, red!75] (9) to[out=135, in=225, looseness=8] (9);
		\draw[thick, red!75, ->, >=stealth] (10) to ++(0.4,+0.4);
		\draw[thick, red!75, ->, >=stealth] (10) to ++(0.4,-0.4);
		
		\node[node, fill=red!50] (11) at (0,-4) {};
		\node[node, fill=red!50] (12) at (1.5,-4) {};
		\draw[thick, red!75] (11) to (12);
		\draw[thick, red!75, ->, >=stealth] (11) to ++(-0.4,+0.4);
		\draw[thick, red!75, ->, >=stealth] (11) to ++(-0.4,-0.4);
		\draw[thick, red!75] (12) to[out=45, in=315, looseness=8] (12);
		
		\node[node] (13) at (0,-5) {};
		\node[node] (14) at (1.5,-5) {};
		\draw[thick] (13) to (14);
		\draw[thick, ->, >=stealth] (13) to ++(-0.4,+0.4);
		\draw[thick, ->, >=stealth] (13) to ++(-0.4,-0.4);
		\draw[thick, ->, >=stealth] (14) to ++(0.4,+0.4);
		\draw[thick, ->, >=stealth] (14) to ++(0.4,-0.4);

		\node[node] (15) at (3.75,-2) {};
		\node[node] (16) at (5.25,-2) {};
		\draw[thick] (15) to[out=30, in=150] (16);
		\draw[thick] (16) to[out=-150, in=-30] (15);
		\draw[thick, ->, >=stealth] (15) to ++(-0.6,0);
		\draw[thick, ->, >=stealth] (16) to ++(0.6,0);

		\node[node] (17) at (7.5,-2) {};
		\node[node] (18) at (9,-2) {};
		\draw[thick] (17) to[out=30, in=150] (18);
		\draw[thick] (17) to (18);
		\draw[thick] (18) to[out=-150, in=-30] (17);
		
	\end{tikzpicture}
	\caption{\label{fig:multigraphplus}The resulting $3$-regular graphs with parallel edges, loops and semi-edges obtained through \texttt{multigraph+}. Only degree distributions yielding graphs are shown. The graphs marked in red are isomorphic.}
\end{figure}

Given a degree distribution $(v_1, \ldots, v_k)$, \texttt{multigraph} will generate all connected graphs with potentially parallel edges that have this degree distribution, i.e., $v_i$ vertices have degree $i$ for each $i \in \{1,\ldots,k\}$. Note that these graphs do not have loops or semi-edges. It is possible to extend \texttt{multigraph} to generate these by, given a regularity $k$ and order $n$, trying all degree distributions $(v_1, \ldots, v_k)$ with $\sum_{i=1}^{n} v_i = n$. The resulting graphs are then completed by adding every possible combination of loops and semi-edges until every vertex has degree $k$. This operation may introduce isomorphic graphs that can be removed using the software package Nauty~\cite{nauty}. An example is shown in Fig.~\ref{fig:multigraphplus}. This extended version of \texttt{multigraph} will be described as \texttt{multigraph+} moving forward.

\subsection{Selection of voltage assignments} \label{subsec:selectionvoltage}
There are multiple ways to exclude specific voltage assignments that will not result in an interesting lift, i.e., a lift satisfying the properties that the different problems are looking for. These exclusions will be broadly categorized as either structure- or assignment-based. Structure-based exclusions only make use of properties of the base graph and group and are completely independent of the chosen voltage assignment. Assignment-based exclusions do make use of the (potentially partial) voltage assignment.

\subsubsection{Structure-based exclusions} \label{subsubsec:structurebasedexclusion}
Since structure-based exclusions do not depend on a voltage assignment, these can be calculated in advance.

\paragraph*{Semi-edge constraints}
For the lift to be interpretable as an undirected graph, the voltage for each dart corresponding with a semi-edge must be its own inverse. Moreover, these darts should not be assigned $0_\Gamma$ as this would result in a loop in the lift. This means a dart corresponding with a semi-edge should only be assigned group elements of order 2. If the chosen group $\Gamma$ has $\ell$ such elements, then each vertex of the base graph may have at most $\ell$ incident semi-edges to avoid using the same voltage multiple times. Combinations of base graphs and groups not adhering to this constraint can be safely skipped as any lift would contain parallel edges.

\paragraph*{Spanning tree and cycles based on the spanning tree}
As stated in Proposition~\ref{prop:tree}, without loss of generality we may pick and fix any spanning tree $T$ and assign $0_\Gamma$ to the darts belonging to this tree. 

This can now be combined with Proposition~\ref{prop:walkgirth} to exclude specific voltages for the remaining darts not in $T$. Consider one such dart $d$ between two distinct vertices $u$ and $v$. Since $T$ is a spanning tree, there must be a path of some length $q$ in $T$ from $u$ to $v$. If this chosen dart is assigned voltage $a$, then the concatenation of the dart $d$ and the path in $T$ from $v$ to $u$ forms a closed, non-reversing walk of length $q+1$ with net voltage $a$. If $a^s = 0_\Gamma$, then the concatenation of this walk with itself $s$ times yields a closed, non-reversing walk of length $(q+1)s$ with net voltage $0_\Gamma$. As this would cause the existence of a cycle of length $(q+1)s$ in the lift, specific voltages for this dart can be excluded based on the girth lower bound to which the generated graphs should adhere.

The above strategy can also be applied to darts belonging to loops. In that case, the concatenation of the loop dart with itself $s$ times results in a closed, non-reversing walk of length $s$ with net voltage $0_\Gamma$, and therefore a cycle of length $s$ in the lift.

If any of the strategies above show that a dart has no voltage that can yield a lift achieving the girth lower bound, then that entire combination of base graph and group can be skipped.

\subsubsection{Assignment-based exclusions}\label{subsubsec:assignmentbasedexclusion}
To facilitate reasoning about which voltage assignments $\alpha$ yield a desirable lift from any given base graph $G$ and group $\Gamma$, it is useful to allow $\alpha$ to be a \emph{partial voltage assignment}. A voltage assignment is partial if it does not yet assign a voltage to all darts. Additionally, define $L_u(G)$ as an ordered \emph{\textbf{l}ist of \textbf{u}seful darts} in $G$. A dart $d$ is considered useful if it is not part of the chosen spanning tree and if $d$ corresponds to a non-semi-edge, then exactly one of $d$ or $d^{-1}$ is considered useful (we can choose this arbitrarily; if $\alpha(d)$ is fixed, then $\alpha(d^{-1})$ is forced since $\alpha(d^{-1})=\alpha(d)^{-1}$).

\paragraph*{Group and edge automorphisms}
It is now possible to apply a voltage assignment $\alpha$ to each dart in $L_u(G)$, resulting in a list of voltages.  Observations~\ref{ob:groupiso}~and~\ref{ob:graphiso} now induce an equivalence relation on all possible lists $\alpha(L_u(G))$ for any given graph $G$. More formally
\[
\alpha(L_u(G)) \sim \alpha'(L_u(G)) \iff \exists \phi_\Gamma, \psi_G: \alpha(L_u(G)) = \phi_\Gamma(\alpha'(\psi_G(L_u(G)))),
\]
with $\phi_\Gamma$ a group automorphism of $\Gamma$, $\phi_G$ an edge automorphism of $G$ that maps every edge outside the chosen spanning tree to another edge outside the chosen spanning tree and $\psi_G$ the natural bijection of darts induced by $\phi_G$. Lists and their corresponding voltage assignments that belong to the same class result in isomorphic lifts. It consequently suffices to keep and construct a lift using a single canonical representative from each equivalence class of assignments, for instance the lexicographically smallest list $\alpha(L_u(G))$. If a partial assignment is rejected as non-canonical, then every completion of it can also be rejected. A formal description of this canonicity verification can be found in Algorithm~\ref{alg:NotCanonical}.

Group and edge automorphisms can be computed using GAP and Nauty~\cite{nauty}, respectively.

\begin{algorithm}[tb]
\caption{Check if the current voltage assignment is canonical under group and edge automorphisms.}
\label{alg:NotCanonical}
\begin{algorithmic}
\Require Graph $G$, group $\Gamma$, partial assignment $\alpha$, ordered dart list $L_u$
\Function{NotCanonical}{$G, \Gamma, \alpha, L_u$}
    \State $\text{VoltageList} \gets \alpha(L_u)$
    \ForAll{edge automorphisms $\phi_G$}
        \State $\psi(G) \gets$ bijection of darts induced by $\phi_G$
        \State $\text{VoltageList}^{\psi_G} \gets$ permute voltages of $\text{VoltageList}$ according to $\psi_G$
        \If{\Call{IsLexSmallerUnderGroup}{$\text{VoltageList}^{\psi_G}, \text{VoltageList}, \Gamma$}}
            \State \Return \textbf{true}
        \EndIf
    \EndFor
    \State \Return \textbf{false}
\EndFunction
\\
\Function{IsLexSmallerUnderGroup}{$a,b,\Gamma$}
    \ForAll{group automorphisms $\phi_\Gamma \in \mathrm{Aut}(\Gamma)$}
        \For{$i=0$ \textbf{to} $|a|-1$}
            \If{$a[i] < \phi_\Gamma [b[i]]$} 
                \State \Return \textbf{true}
            \ElsIf{$a[i] > \phi_\Gamma [b[i]]$} \State \textbf{break}
            \EndIf
        \EndFor
    \EndFor
    \State \Return \textbf{false}
\EndFunction
\end{algorithmic}
\end{algorithm}

\paragraph*{Cycles based on partial voltage assignments}
Similar to the structure-based exclusions, Proposition~\ref{prop:walkgirth} can be used to verify whether $\alpha$ results in a lift exhibiting the desired minimum girth $\gmin$. It suffices to generate all closed non-reversing walks with a length strictly shorter than $\gmin$. If $\alpha$ is partial, then the darts that did not yet receive a voltage should not be used in these walks. If any of these walks have a net voltage of $0_\Gamma$ then this voltage assignment can be skipped. Again, if a partial voltage assignment $\alpha$ is excluded, then any completion of $\alpha$ will necessarily contain the same short closed walk and can therefore also be rejected.

The verification of potentially short walks can be made more efficient by exploiting the incremental nature of assigning voltages. Suppose $\alpha'$ is identical to $\alpha$ except for the assignment on a single dart $d$ and $\alpha$ was already checked for short walks. Since all walks not using $d$ were already verified under $\alpha$, it suffices to check only those walks that include $d$. Moreover, the starting dart and direction of a walk respectively yield a conjugate and inverse group element as net voltage. They consequently do not influence the order of the resulting net voltage \cite[p.63]{topologicalGraphTheory}, which means only walks that start with $d$ are of interest. This efficient version is formally described in Algorithm~\ref{alg:CannotAchieveGirth}.

In the special case where $\alpha$ assigns no voltages yet to darts in $L_u(G)$, the procedure reduces to the structure-based exclusion based on the spanning tree.

\begin{algorithm}[tb]
\caption{Check if the current partial assignment can still achieve the minimum girth by generating all closed non-reversing walks from a given dart.}
\label{alg:CannotAchieveGirth}
\begin{algorithmic}
\Require Graph $G$, group $\Gamma$, dart $d=(u,v)$ with the most recent assignment, minimum girth $g_{\min}$, and a potentially partial assignment $\alpha$
\Function{CannotAchieveGirth}{$G, \Gamma, d, g_{\min}, \alpha$}
    \State \Return \Call{ExistsCNRWalk}{$\alpha(d^{-1}), d^{-1}, u, v, 1, \gmin-1, \alpha, G, \Gamma$}
\EndFunction
\\
\Function{ExistsCNRWalk}{$a_{\text{cur}}, d_{\text{prev}}, v_{\text{cur}}, v_{\text{target}}, \ell, \ell_{\max}, \alpha, G, \Gamma$}
    \If{$\ell \geq \ell_{\max}$}
        \State \Return \textbf{false}
    \EndIf
    \ForAll{darts $d=(v_{\text{cur}},w)$ in $G$}
        \If{$\alpha(d)$ is not defined or $d = d_{\text{prev}}^{-1}$}
            \State \textbf{continue}
        \EndIf
        \State $a_{\text{new}} \gets a_{\text{cur}} \cdot \alpha(d)$
        \If{$w = v_{\text{target}}$ \textbf{and} $a_{\text{new}} = 0_\Gamma$}
            \State \Return \textbf{true}
        \EndIf
        \If{\Call{ExistsCNRWalk}{$a_{\text{new}}, d, w, v_{\text{target}}, \ell+1, \ell_{\max}, \alpha, G, \Gamma$}}
            \State \Return \textbf{true}
        \EndIf
    \EndFor
    \State \Return \textbf{false}
\EndFunction
\end{algorithmic}
\end{algorithm}

\subsection{Efficient filtering}\label{subsec:efficientfiltering}
Only connected lifts that are non-isomorphic to previously generated graphs and that satisfy the required girth (and optionally edge- or vertex-girth-regularity conditions) are of interest. These properties can be checked efficiently after constructing the lift. More specifically, a canonical representation of the lift can be constructed using Nauty to filter out isomorphic copies, while both the girth and edge- or vertex-girth-regularity conditions can be calculated using the algorithm described in \cite{goedgebeur2024egrgraphs}.

\subsection{Iterating over voltage assignments}\label{subsec:iteratingvoltages}
Once the admissible base graph-group combinations have been determined and the structure- and assignment-based exclusions have been applied, the remaining search space of voltage assignments can be explored.

\subsubsection{Exhaustive backtracking} \label{sec:BTA}
Algorithm~\ref{alg:BTA} describes a BackTracking Algorithm (BTA) that recursively assigns voltages to darts in $L_u(G)$. At each step it applies the exclusions from Subsections~\ref{subsubsec:structurebasedexclusion} and \ref{subsubsec:assignmentbasedexclusion} to prune infeasible or non-canonical assignments. Every surviving assignment is lifted to a graph, which is then filtered using the methods of Section~\ref{subsec:efficientfiltering}.

\begin{algorithm}[tb]
\caption{A recursive backtracking algorithm for generating lifts achieving a given minimum girth and using canonical voltage assignments.}
\label{alg:BTA}
\begin{algorithmic}
\Require A graph $G$, a group $\Gamma$, and a minimum girth $g_{\text{min}}$
\Function{BTA}{$G, \Gamma, g_{\text{min}}$}
    \State $(N, L_u, \alpha) \gets$ \Call{StructuralChecks}{$G, \Gamma, g_{\text{min}}$} \Comment{Abort if impossible}
    \State \Call{BTARecursive}{$0, |L_u|, N, L_u, \alpha, G, \Gamma, g_{\text{min}}$}
\EndFunction
\\
\Function{BTARecursive}{$\ell, \ell_{\max}, N, L_u, \alpha, G, \Gamma, g_{\text{min}}$}
    \If{$\ell = \ell_{\max}$}
        \State $G^\alpha \gets$ \Call{Lift}{$G, \Gamma, \alpha$}
        \If{\Call{IsInteresting}{$G^\alpha$}}
            \State \Call{Output}{$G^\alpha$}
        \EndIf
        \State \Return
    \EndIf

    \State $d \gets L_u[\ell]$
    \ForAll{voltages $a \in N(d)$}
        \State $\alpha(d) \gets a,\;\; \alpha(d^{-1}) \gets a^{-1}$

        \If{\Call{CannotAchieveGirth}{$G, \Gamma, d, g_{\text{min}}, \alpha$} 
            \textbf{or} \Call{NotCanonical}{$G, \Gamma, \alpha, L_u$}}
            \State \textbf{continue}
        \EndIf

        \State \Call{BTARecursive}{$\ell+1, \ell_{\max}, N, L_u, \alpha, G, \Gamma, g_{\text{min}}$}
    \EndFor
    \State remove $\alpha(d), \alpha(d^{-1})$ \Comment{Algorithm will backtrack, so restore initial state}
\EndFunction
\end{algorithmic}
\end{algorithm}

\subsubsection{Heuristic search} \label{sec:TS}
While BTA is exhaustive and guarantees completeness, it can become computationally expensive as the orders of the group and base graph grow. As an alternative, metaheuristics such as Tabu Search (TS)~\cite{glover1989} can be used to sample promising voltage assignments without fully exploring the search space. TS starts from an initial assignment and iteratively explores neighboring assignments, guided by a cost function reflecting the desired properties (e.g.\ large girth, vertex-girth-regularity or edge-girth-regularity). The tabu list and perturbations help avoid revisiting solutions and escaping local optima. We now explain the different elements used in this algorithm and refer the reader to Algorithm~\ref{alg:TS} for the pseudocode.

\paragraph*{Neighboring assignments}
The neighbors of a voltage assignment $\alpha$ consist of all alternative assignments that differ from $\alpha$ in the group element assigned to the darts corresponding with exactly one edge. An initial assignment can be obtained with BTA by simply stopping once the first valid voltage assignment is found.

A significant portion of these candidates can be excluded based on assignment-dependent strategies as described in Subsection~\ref{subsubsec:assignmentbasedexclusion}. For example, one may again require that assignments are canonical according to Algorithm~\ref{alg:NotCanonical} and that they achieve a certain minimum girth via Algorithm~\ref{alg:CannotAchieveGirth}. Note that, for the cage problem specifically, the minimum girth $g_{\text{nbr}}$ when selecting candidate neighbors should be less than the target minimum girth of the lift $\gmin$. Picking a value for $g_{\text{nbr}}$ close to $\gmin$ prunes more bad candidates, but also risks pruning all available candidates and makes the search for an initial assignment take longer.

\paragraph*{Cost functions}
Two cost functions are considered in this work. The first cost function is based on the girth and utilizes the intuition that short closed walks of a low net voltage order are undesirable. For a walk of length $q$ with net voltage $a$, the cost contribution is defined as
\[
    f(q, a)=
            \begin{cases}
                \frac{1}{q \cdot \Ord(a)} &\text{if }\Ord(a) \neq 1,\\[4pt]
                C &\text{otherwise},
            \end{cases}
\]
with $C$ a large penalty constant for walks with net voltage $0_\Gamma$. The algorithm will not exhaustively generate all walks up to a certain length, but will instead sample $m$ walks. The total cost of an assignment based on $m$ sampled walks is
\[
    c_{\text{girth}} = \sum_{i=1}^m f(q_i, a_i).
\]

The second cost function is used when trying to find small $vgr$- and $egr$-graphs. Since the target minimum girth $\gmin$ is usually smaller compared to the cage problem, many voltage assignments will achieve this $\gmin$. Consequently, most time is spent repeatedly verifying whether or not the resulting lifts are  $vgr$- and $egr$-graphs. To approximate this, one can compute how often each vertex or dart appears in closed non-reversing walks of net voltage $0_\Gamma$ with length equal to the girth, and these frequencies should be uniform. Let the resulting frequencies be $f_V$ (for vertices) and $f_D$ (for darts) then a suitable cost function is the minimum of their variances:
\[
    c_{\text{reg}} =\min\Bigg[\frac{1}{|f_V|}\sum_{i=1}^{|f_V|} (f_{V,i} - \bar{f_V})^2,\quad \frac{1}{|f_D|}\sum_{i=1}^{|f_D|} (f_{D,i} - \bar{f_D})^2\Bigg].
\]

\paragraph*{Tabu list}
To avoid getting stuck in local optima, TS employs a so-called \emph{tabu list}. This list prevents the algorithm from immediately returning to previously visited configurations, thereby encouraging exploration of diverse voltage assignments. In this context the tabu list stores modifications rather than entire assignments: specifically,
we store which modification would again undo the most recent modification.

Since neighbors are obtained by changing exactly one voltage, a neighbor can be skipped if its modification is marked as tabu, i.e., occurs in the tabu list. An exception is made if this neighbor improves upon the globally best assignment found so far. In that case, the move is allowed even if it occurs in the tabu list. The tabu list has a fixed length, and once full, newly added moves overwrite the oldest entries, ensuring that prohibitions are only temporary and eventually expire.

\paragraph*{Perturbations}
It may still happen that TS cannot escape a local optimum, or that all possible modifications are in the tabu list. In such cases it is useful to perturb the current solution, for example by applying several random modifications. Whether the algorithm is stuck in a local optimum can be judged from the number of iterations since the last improvement of the global best solution. Perturbations hopefully push TS into a previously unexplored region of the search space. 

\begin{algorithm}[tb]
\caption{Constructs a subset of all lifts based on TS.}\label{alg:TS}
\begin{algorithmic}
\Require A graph $G$, a group $\Gamma$, a minimum girth for the lift $g_{\text{min}}$ and a minimum girth for the neighbor solutions $g_{\text{nbr}}$
\Function{TS}{$G, \Gamma, \gmin, g_{\text{nbr}}$}
\State $N, L_u, \alpha \gets $\Call{StructureBasedChecks}{$G, \Gamma, g_{\text{min}}$} \Comment{Abort if impossible}
\State Use \Call{BTARecursive}{$0, |L_u|, N, L_u, \alpha, G, \Gamma, g_{\text{nbr}}$} for an initial solution

\State $L_{\text{tabu}} \gets \emptyset$

\State $c_{\text{global}} \gets \infty$
\While{no stopping criterion is satisfied}
    \If{!\Call{CannotAchieveGirth}{$G, \Gamma, d, \gmin, \alpha$}}
        \State $G^\alpha \gets$ \Call{Lift}{$G, \Gamma, \alpha$} 
        \If{\Call{IsInteresting}{$G^\alpha$}}
            \State \Call{Output}{$G^\alpha$}
        \EndIf
    \EndIf

    \If{too many iterations without improvement of $c_{\text{global}}$}
        \State Perturb
    \EndIf
    
    \State $c_{\text{best}} \gets \infty$
    \ForAll{darts $d\in L_u$}
        \State $a \gets \alpha(d)$
        \ForAll{voltages $b \in N(d)$}
            \If{\Call{CannotAchieveGirth}{$G, \Gamma, d, g_{\text{nbr}}, \alpha$}
                \Statex \hspace{\algorithmicindent} \textbf{or} \Call{NotCanonical}{$G, \Gamma, \alpha, L_u$}
            }
                \State \textbf{Continue}
            \EndIf
            \State $\alpha(d) \gets b, \alpha(d^{-1}) \gets b^{-1}$
            \State $c \gets $ \Call{CostFunction}{$G, \Gamma, \alpha$}
            \If{$c < c_{\text{best}}$ \textbf{and} ($(d, b) \notin L_{\text{tabu}}$ \textbf{or} $c < c_{\text{global}}$)}
                \State $c_{\text{best}} \gets c, d_{\text{best}} \gets d, b_{\text{best}} \gets b$
            \EndIf
        \EndFor
        \State $\alpha(d) \gets a, \alpha(d^{-1}) \gets a^{-1}$
    \EndFor

    \State $c_{\text{global}} \gets \min(c_{\text{global}}, c_{\text{best}})$
    \If{$c_{\text{best}} \neq \infty$}
        \State $\alpha(d_{\text{best}}) \gets b_{\text{best}}, \alpha(d_{\text{best}}^{-1}) \gets b_{\text{best}}^{-1}$
        \State $L_{\text{tabu}} \gets L_{\text{tabu}} \cup \{(d, b_{\text{best}}^{-1})\}$
        \State Expired entries are dropped from $L_{\text{tabu}}$
    \Else 
        \State Perturb
    \EndIf    
\EndWhile
\EndFunction
\end{algorithmic}
\end{algorithm}

\subsubsection{Hyperparameters}\label{subsec:hyperparam}
Before running experiments, suitable values must be chosen for the hyperparameters that guide the search process (e.g., time limits and cost-function weights). We distinguish between parameters relevant for $vgr$- and $egr$-graphs and those for the remaining problems, denoted respectively by $\mathcal{H}_{\text{reg}}$ and $\mathcal{H}_{\text{girth}}$. An overview of the chosen values is provided in Table~\ref{tab:hyperparam}.

For both edge and group automorphisms, the cutoffs were determined empirically based on the resulting time savings. Considering too many automorphisms can increase runtime, as checking them may take longer than constructing redundant lifts. Analogously, only limited time is allocated for detecting useful graph automorphisms.

In all experiments, BTA and TS were given fixed per-instance time budgets. TS was only invoked when BTA was not able to explore the whole search space within the given time limit.

\begin{table}
    \centering
    \caption{Chosen hyperparameters. Time limits are expressed in seconds.}
    \label{tab:hyperparam}
    \begin{tabular}{p{2.8cm} p{7.2cm} p{1.2cm} p{1.2cm}}
    \toprule
    Hyperparameter & Description & $\mathcal{H}_{\text{girth}}$  & $\mathcal{H}_{\text{reg}}$\\ \midrule
         $\#\phi_G$ & Maximum number of edge automorphisms. & 200 & 200 \\
         $\#\phi_\Gamma$ & Maximum number of group automorphisms. & 2000 & 2000 \\ \midrule
         t$_{\phi_G}$ & Time limit for searching useful edge automorphisms. & 2 & 2 \\
         t$_{\text{BTA}}$ & Time allocated to BTA per graph–group combination. & 20 & 20 \\
         t$_{\text{TS}}$ & Time allocated to TS per graph–group combination. & 20 & 20 \\
         t$_{\text{TS, init}}$ & Time limit for TS to find an initial solution. & 2 & 2 \\\midrule
         $g_{\text{nbr}}$ & Minimum girth of neighbor solutions in TS. & $g_{\text{min}}$-2 & $g_{\text{min}}$\\
         $|L_{\text{tabu}}|$ & Size of the tabu list. & $3|\Gamma|$ & $3|\Gamma|$ \\
         Cost$_{\text{TS}}$ & Cost function used by TS. & $c_{\text{girth}}$ & $c_{\text{reg}}$ \\
         $C$ & Penalty constant for short walks with net voltage $0_\Gamma$ in $c_{\text{girth}}$. & 1000 & / \\
         $\#W$ & Number of walks sampled for $c_{\text{girth}}$. & 500 & / \\
         \bottomrule
    \end{tabular}
\end{table}

\subsection{Hill climbing heuristic} \label{subsec:hill climbing}

The algorithms discussed so far fix the target girth $g$, a given base graph and group and then search for a suitable voltage assignment. The hill climbing heuristic discussed in the current subsection takes a different approach: it fixes the target girth $g$, the order $n$ of the base graph (but not the base graph itself) and the group and then simultaneously searches for a suitable base graph of order $n$ and a voltage assignment using the given group. 

More precisely, the hill climbing heuristic starts from a graph consisting of $n$ isolated vertices. This graph will serve as the base graph and will be modified in each iteration of the algorithm by adding a pair of opposing darts and group elements assigned to these darts such that the resulting lift has girth at least $g$. Let $G$ be the current graph, let $\mathcal{T}(G)$ be the set of all tuples consisting of opposing darts and group elements assigned to these opposing darts such that the lift of the graph obtained by adding these darts and voltage assignment has girth at least $g$. The hill climbing heuristic tries to maximize the number of opposing darts that will eventually be added. Therefore, the algorithm will consider for each tuple $t \in \mathcal{T}(G)$ the graph $G_t$ obtained by modifying $G$ using the tuple $t$ (i.e., adding the corresponding darts and voltage assignment) and greedily chooses the tuple $t$ that maximizes $|\mathcal{T}(G_t)|$ (i.e., modifications that lead to a graph where many new modifications are possible). After each iteration, the hill climbing heuristic considers the lift of the current graph and outputs this graph in case it breaks any of the records. We remark that this hill climbing heuristic can again apply most of the optimizations that were previously discussed (e.g.\ for efficiently computing the resulting girth or ruling out voltage assignments).

\subsection{Excision}\label{subsec:excision}

Excision is the technique of removing a subset of vertices from a small known $(k,g)$-graph and adding edges to the resulting graph to obtain a $k$-regular graph of large girth (typically $g-1$). Notably, a lot of the current record graphs were obtained by using this technique \cite{cagesurvey}. We obtain six improved upper bounds using excision techniques:
\[
n(4,11) \leq 713, \quad 
n(6,11) \leq 7783,\quad 
n(8,7) \leq 774, \quad 
n(10,7) \leq 1608,\quad 
n(12,7) \leq 2890, \quad 
n(14,7) \leq 4716.
\]
We obtained these bounds by first deleting a specific set of vertices and its incident edges from a small known $(k,g+1)$-graph. Let us call the obtained non-regular graph $G'$. Secondly, we feed $G'$ to a $(k,g)$-graph generation algorithm, as described in \cite{EXOO2011166,10.5555/314613.314699}, which adds edges to $G'$ in order to complete it to a $k$-regular graph of girth at least $g$ with the same number of vertices (when completion is possible).
Which vertices we remove from the starting graph, depends on the different cases. Let us denote the set of vertices at distance $i$ from vertex $v$ as $N_i(v)$, and the set of vertices at distance $i$ from both vertex $u$ and $v$ as $N_i(u,v)= N_i(u) \cap N_i(v)$. 

For girth 8, we slightly modified the excised set from the one described for even degrees $4 \le k \le 14$ by de~Ruiter and Biggs~\cite{deruiter2015applications}. They remove two vertices $u$ and $v$ at distance 4, their neighborhoods $N_1(u)$, $N_1(v)$ and $N_2(u,v)$ -- except three vertices of $N_2(u,v)$. Since all $(k,8)$-cages with even degree $4 \le k \le 14$ have the property that $|N_2(u,v)|=k$ \cite{deruiter2015applications}, it follows that the excised set has size $3k-1$. For even degrees $8 \le k \le 14$, however, removing the three vertices from the excised set is unnecessary. When fully excising $N_2(u,v)$ (in addition to $u$, $v$, $N_1(u)$ and $N_1(v)$), the $(k,g)$-graph generator was able to complete the graph to a $k$-regular graph of girth 7. Hence, we improved the upper bound on $n(k,7)$ for even $8 \le k \le 14$ by three vertices.

\begin{figure}
\centering
\begin{tikzpicture}[
  every node/.style={circle,draw, minimum size=0.4cm},
  edge/.style={thick,black}
]

\def\xsep{1.0}
\def\ysep{1.0}

\coordinate (row1) at (0,6*\ysep);  
\coordinate (row2) at (0,5*\ysep);
\coordinate (row3) at (0,4*\ysep);
\coordinate (row4) at (0,3*\ysep);
\coordinate (row5) at (0,2*\ysep);
\coordinate (row6) at (0,1*\ysep); 
\coordinate (row7) at (0,0*\ysep); 

\node[fill=red!50] (T) at (row1) {};
\node[fill=red!50] (B) at (row7) {};

\foreach \r/\y in {2/5, 3/4, 4/3, 5/2, 6/1} {
  \foreach \c in {1,2,3,4} {
   
    \pgfmathsetmacro{\x}{(\c-2.5)*\xsep}
    \ifthenelse{\r = 2 \OR \r=4 \OR \r=6 \OR \(\r=3 \AND \c=1\)}{\node[fill=red!50] (v\y\c) at (\x,\y*\ysep) {};}{\node (v\y\c) at (\x,\y*\ysep) {};}
  }
}

\foreach \c in {1,2,3,4}{
  \draw[edge] (T) -- (v5\c);
}

\foreach \c in {1,2,3,4}{
  \draw[edge] (B) -- (v1\c);
}

\foreach \c in {1,2,3,4}{
  \draw[edge] (v5\c) -- (v4\c);
  \draw[edge] (v4\c) -- (v3\c);
  \draw[edge] (v3\c) -- (v2\c);
  \draw[edge] (v2\c) -- (v1\c);
}

\foreach \c in {1,2,3,4}{
  \draw[edge] (v5\c) -- (v4\c);
  \draw[edge] (v4\c) -- (v3\c);
  \draw[edge] (v3\c) -- (v2\c);
  \draw[edge] (v2\c) -- (v1\c);
}

\node[draw=none] (u) at (4,6*\ysep) {$u$};
\node[draw=none] (N1u) at (4,5*\ysep) {$N_1(u)$};
\node[draw=none] (N2u) at (4,4*\ysep) {$N_2(u) \cap N_4(v)$};
\node[draw=none] (N3) at (4,3*\ysep) {$N_3(u,v)$};
\node[draw=none] (N2v) at (4,2*\ysep) {$N_4(u) \cap N_2(v)$};
\node[draw=none] (N1v) at (4,1*\ysep) {$N_1(v)$};
\node[draw=none] (u) at (4,0*\ysep) {$v$};

\node[fill=red!50] (T') at (9,6*\ysep) {};
\node[fill=red!50] (B') at (9,0*\ysep) {};

\foreach \r/\y in {2/5, 3/4, 4/3, 5/2, 6/1} {
  \foreach \c in {1,2,3,4,5,6} {
    \pgfmathsetmacro{\x}{9+(\c-3.5)*\xsep}
    \ifthenelse{\r=4 \AND \c<4}{\node (v'\y\c) at (\x,\y*\ysep) {};}{\node[fill=red!50] (v'\y\c) at (\x,\y*\ysep) {};}
  }
}

\foreach \c in {1,2,3,4,5,6}{
  \draw[edge] (T') -- (v'5\c);
}

\foreach \c in {1,2,3,4,5,6}{
  \draw[edge] (B') -- (v'1\c);
}

\foreach \c in {1,2,3,4,5,6}{
  \draw[edge] (v'5\c) -- (v'4\c);
  \draw[edge] (v'4\c) -- (v'3\c);
  \draw[edge] (v'3\c) -- (v'2\c);
  \draw[edge] (v'2\c) -- (v'1\c);
}

\foreach \c in {1,2,3,4,5,6}{
  \draw[edge] (v'5\c) -- (v'4\c);
  \draw[edge] (v'4\c) -- (v'3\c);
  \draw[edge] (v'3\c) -- (v'2\c);
  \draw[edge] (v'2\c) -- (v'1\c);
}

\end{tikzpicture}
\caption{The excised set of vertices indicated in red from the $(4,12)$-cage (left) and the $(6,12)$-cage (right).}
\label{fig:excsionGirth12}
\end{figure}
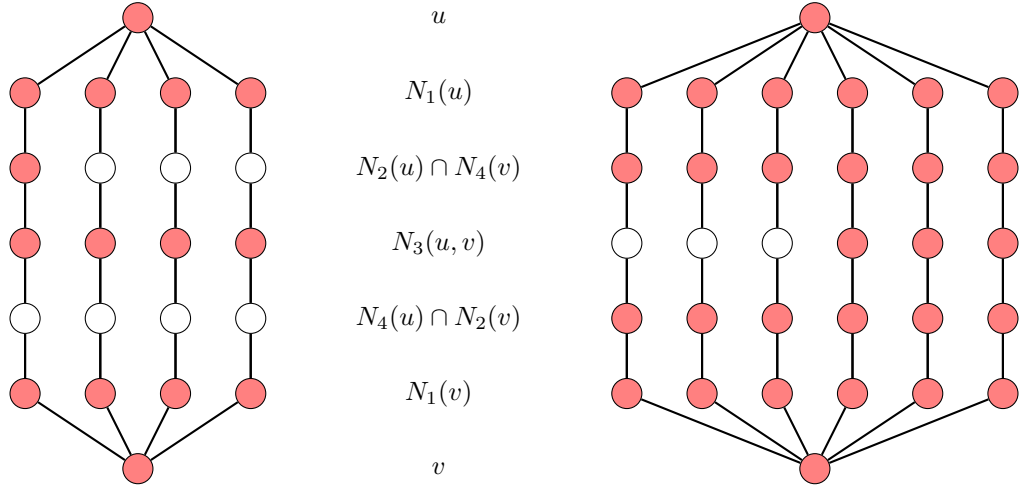

Inspired by de~Ruiter and Biggs' work, we devised excision sets for the known $(4,12)$- and $(6,12)$-cage, depicted in Fig.~\ref{fig:excsionGirth12}. First, we note that for two vertices $u$, $v$ at distance 6, we have $|N_3(u,v)|=|N_2(u) \cap N_4(v)|=|N_4(u) \cap N_2(v)|=k$. For the $(4,12)$-cage, we remove the two vertices $u$ and $v$ at distance 6, their neighborhoods $N_1(u)$, $N_1(v)$, the complete set $N_3(u,v)$ and one vertex from either $N_2(u) \cap N_4(v)$ or $N_4(u) \cap N_2(v)$ (both sets work because of the symmetry of the $(4,12)$-cage). We thus excise $3k+3$ vertices. For $(6,12)$, we can follow a similar approach. However, we can excise an even larger set of $5k-1$ vertices. This set consists of vertices $u$ and $v$ at distance 6, $N_1(u)$, $N_1(v)$, $N_2(u) \cap N_4(v)$, $N_4(u) \cap N_2(v)$ and $N_3(u,v)$ -- except three vertices from $N_3(u,v)$.

\section{Results}\label{sec:results}
The total computational effort underlying this work corresponds to roughly five CPU-years, carried out on the infrastructure of the Flemish Supercomputer Center. All code and data involving lifts or excision are publicly accessible at \url{https://github.com/AGT-Kulak/smallRegGirthGraphs}. 

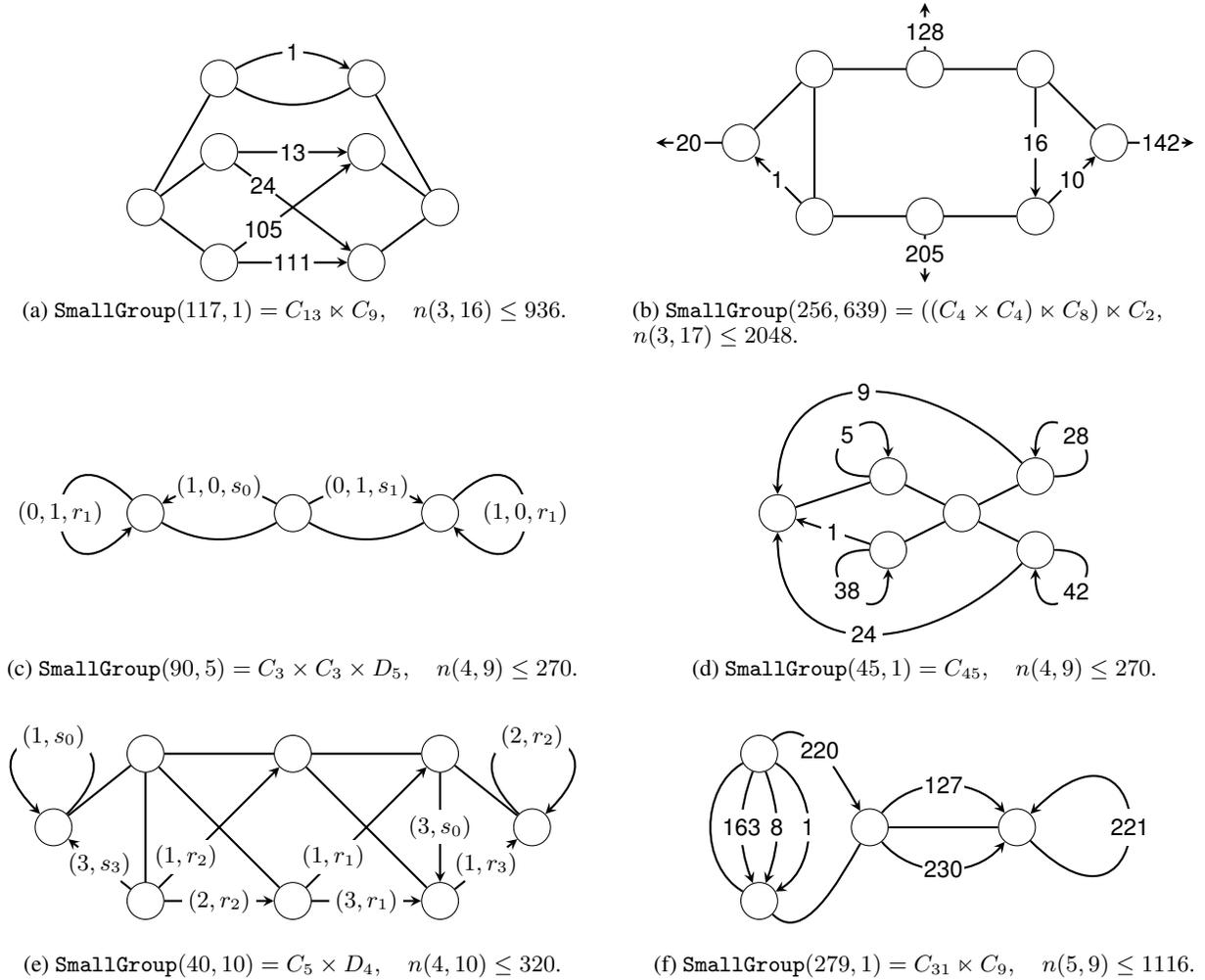
\begin{figure}
\centering
\begin{subfigure}{.48\linewidth}
    \centering
    \begin{tikzpicture}[
            node/.style={circle, draw, minimum size=0.5cm},
            label/.style={midway, fill=white, font=\sffamily\small}
        ]
        
        \node[node] (A) at (-1,0.75) {};
        \node[node] (B) at (1,0.75) {};
        
        \node[node] (C) at (-1,-0.75) {};
        \node[node] (D) at (1,-0.75) {};
        
        \node[node] (E) at (-2,0) {};
        \node[node] (F) at (2,0) {};

        \node[node] (G) at (1,1.75) {};
        \node[node] (H) at (-1,1.75) {};
        
        \draw[thick, ->, >=stealth] (A) -- node[label, circle, inner sep=0pt] {13} (B);
        \draw[thick, ->, >=stealth] (A) -- node[label, xshift=-0.4cm, yshift=0.30cm, circle, inner sep=0pt] {24} (D);
        \draw[thick, ->, >=stealth] (C) -- node[label, xshift=-0.4cm, yshift=-0.30cm, circle, inner sep=0pt] {105} (B);
        \draw[thick, ->, >=stealth] (C) -- node[label, circle, inner sep=0pt] {111} (D);

        \draw[thick] (E) --  (A);
        \draw[thick] (E) --  (C);
        \draw[thick] (F) --  (B);
        \draw[thick] (F) --  (D);

        \draw[thick] (E) --  (H);
        \draw[thick] (F) --  (G);

        \draw[thick, ->, >=stealth] (H) to[out=30, in=150] node[label, circle, inner sep=0pt] {1} (G);
        \draw[thick] (G) to[out=-150, in=-30] (H);
        \end{tikzpicture}
        \caption{$\texttt{SmallGroup}(117, 1) = C_{13} \ltimes C_9$ and $n(3, 16) \leq 936$.\\\phantom{.}}
        \label{fig:cageRecord316}
\end{subfigure}%
    \hfill
\begin{subfigure}{.48\linewidth}
    \centering
    \begin{tikzpicture}[
            node/.style={circle, draw, minimum size=0.5cm},
            label/.style={midway, fill=white, font=\sffamily\small}
        ]
        
        \node[node] (A) at (-2.5,0) {};
        
        \node[node] (B) at (-1.5,1) {};
        \node[node] (C) at (-1.5,-1) {};

        \node[node] (D) at (0,1) {};
        \node[node] (E) at (0,-1) {};
        
        \node[node] (F) at (1.5,1) {};
        \node[node] (G) at (1.5,-1) {};

        \node[node] (H) at (2.5,0) {};
        
        \draw[thick] (A) --  (B);
        \draw[thick, <-, >=stealth] (A) -- node[label, circle, inner sep=0pt] {1}  (C);
        
        \draw[thick] (B) --  (C);
        
        \draw[thick] (B) --  (D);
        \draw[thick] (C) --  (E);

        \draw[thick] (D) --  (F);
        \draw[thick] (E) --  (G);

        \draw[thick, <-, >=stealth] (G) -- node[label, circle, inner sep=0pt] {16} (F);
        
        \draw[thick] (F) -- (H);
        \draw[thick, ->, >=stealth] (G) -- node[label, circle, inner sep=0pt] {10} (H);

        \draw[thick, ->, >=stealth] (A) to node[label, circle, inner sep=0pt] {20} ++(-1.15, 0);
        \draw[thick, ->, >=stealth] (H) to node[label, circle, inner sep=0pt] {142} ++(1.15, 0);
        \draw[thick, ->, >=stealth] (D) to node[label, circle, inner sep=-2pt, yshift=-2pt] {128} ++(0, 0.9);
        \draw[thick, ->, >=stealth] (E) to node[label, circle, inner sep=-2pt, yshift=2pt] {205} ++(0, -0.9);
    \end{tikzpicture}
    \caption{$\mathtt{SmallGroup}(256,639) = ((C_{4} \times C_4) \ltimes C_8) \ltimes C_2$ and $ n(3, 17) \leq 2048$.}
    \label{fig:cageRecord317}
    \vspace{1em}
\end{subfigure}%

\bigskip
\begin{subfigure}{.48\linewidth}
    \centering
    \begin{tikzpicture}[
            node/.style={circle, draw, minimum size=0.5cm},
            label/.style={midway, fill=white, font=\sffamily\small}
        ]

        \clip(-3.80,-1.75) rectangle (3.80,1.75);

        \begin{scope}
        \node[node] (A) at (-2,0) {};
        \node[node] (C) at (0,0) {};
        \node[node] (B) at (2,0) {};
        
        \draw[<-, thick, >=stealth] (A) to[out=30, in=150] node[label, circle, inner sep=0pt] {$(1, 0, s_0)$} (C);
        \draw[thick] (C) to[out=-150, in=-30] (A);

        \draw[thick, >=stealth, ->] (C) to[out=30, in=150] node[label, circle, inner sep=0pt] {$(0, 1, s_1)$} (B);
        \draw[thick] (B) to[out=-150, in=-30] (C);

        \draw[thick, ->, >=stealth] (A) to[out=135, in=225, looseness=13] node[label, inner sep=0pt, minimum size=0.5cm] {$(0, 1, r_1)$} (A);
        
        \draw[thick, ->, >=stealth] (B) to[out=45, in=315, looseness=13] node[label, inner sep=0pt, minimum size=0.5cm] {$(1, 0, r_1)$} (B);
        \end{scope}
    \end{tikzpicture}
        \caption{$\mathtt{SmallGroup}(90, 5) = C_3 \times C_3 \times D_5$ and $n(4, 9) \leq 270$.}
        \label{fig:cageRecord49a}
\end{subfigure}%
    \hfill
\begin{subfigure}{.48\linewidth}
    \centering
    \begin{tikzpicture}[
            node/.style={circle, draw, minimum size=0.5cm},
            label/.style={midway, fill=white, font=\sffamily\small, transform shape=false}
        ]

        \clip(-3,-1.75) rectangle (2,1.75);
        
        \node[node] (0) at (0,0) {};
        
        \node[node] (1) at (-1,0.5) {};
        \node[node] (2) at (1,0.5) {};
        \node[node] (3) at (1,-0.5) {};
        \node[node] (4) at (-1,-0.5) {};

        \node[node] (5) at (-2.5,0) {};

        \draw[thick] (0) -- (1);
        \draw[thick] (0) -- (2);
        \draw[thick] (0) -- (3);
        \draw[thick] (0) -- (4);

        \draw[thick] (1) -- (5);
        \draw[thick, ->, >=stealth] (4) -- node[label, inner sep=2pt] {1} (5);

        \draw[thick, <-, >=stealth] (3) to[out=-90, in=0, looseness=8] node[label, inner sep=2pt] {42} (3);
        \draw[thick, ->, >=stealth] (2) to[out=0, in=90, looseness=8] node[label, inner sep=2pt] {28} (2);
        \draw[thick, <-, >=stealth] (1) to[out=90, in=180, looseness=8] node[label, inner sep=2pt] {5} (1);
        \draw[thick, ->, >=stealth] (4) to[out=180, in=270, looseness=8] node[label, inner sep=2pt] {38} (4);

        \draw[thick, ->, >=stealth] (2) to[out=135, in=90, looseness=1.4] node[label, inner sep=2pt] {9} (5);
        \draw[thick, ->, >=stealth] (3) to[out=225, in=270, looseness=1.4] node[label, inner sep=2pt] {24} (5);
        
    \end{tikzpicture}
        \caption{$\texttt{SmallGroup}(45, 1)= C_{45}$ and $n(4, 9) \leq 270$.}
        \label{fig:cageRecord49b}
\end{subfigure}%

\bigskip
\begin{subfigure}{.48\linewidth}
    \centering
    \begin{tikzpicture}[
            node/.style={circle, draw, minimum size=0.5cm},
            label/.style={midway, fill=white, font=\sffamily\small}
        ]
        
        \node[node] (4) at (0.75,0) {};
        \node[node] (0) at (2, 1) {};
        \node[node] (1) at (2,-1) {};
        \node[node] (2) at (4, 1) {};
        \node[node] (3) at (4,-1) {};
        \node[node] (5) at (6, 1) {};
        \node[node] (6) at (6,-1) {};
        \node[node] (7) at (7.25,0) {};

        \clip(0.1,1.5) rectangle (7.9,-1.5);
        
        \draw[thick, >=stealth, ->] (4) to[out=45, in=135, looseness=14] node[label, inner sep=2pt] {$(1, s_0)$} (4);
        \draw[thick, >=stealth, <-] (7) to[out=405, in=135, looseness=14] node[label, inner sep=2pt] {$(2, r_2)$} (7);

        \draw[thick] (4) -- (0);
        \draw[thick, >=stealth, <-] (4) -- node[label, circle, inner sep=-2pt] {$(3, s_3)$} (1);
        \draw[thick] (0) -- (1);
        \draw[thick] (0) -- (2);
        \draw[thick] (0) -- (3);
        \draw[thick, >=stealth, ->] (1) -- node[label, inner sep=2pt, xshift=-0.45cm, yshift=-0.4cm] {$(1, r_2)$} (2);
        \draw[thick, >=stealth, ->] (1) -- node[label, inner sep=2pt] {$(2, r_2)$} (3);
        \draw[thick] (2) -- (5);
        \draw[thick] (2) -- (6);
        \draw[thick, >=stealth, ->] (3) -- node[label, inner sep=2pt, xshift=-0.45cm, yshift=-0.4cm] {$(1, r_1)$} (5);
        \draw[thick, >=stealth, ->] (3) -- node[label, inner sep=2pt] {$(3, r_1)$} (6);
        \draw[thick, >=stealth, ->] (5) -- node[label, inner sep=2pt] {$(3, s_0)$} (6);
        \draw[thick] (5) -- (7);
        \draw[thick, >=stealth, ->] (6) -- node[label, circle, inner sep=-2pt] {$(1, r_3)$} (7);
    \end{tikzpicture}
        \caption{$\mathtt{SmallGroup}(40, 10) = C_5 \times D_{4}$ and $n(4, 10) \leq 320$.}
        \label{fig:cageRecord410}
\end{subfigure}%
    \hfill
\begin{subfigure}{.48\linewidth}
    \centering
    \begin{tikzpicture}[
            node/.style={circle, draw, minimum size=0.5cm},
            label/.style={midway, fill=white, font=\sffamily\small, transform shape=false}
        ]

       \node[node] (A) at (-2.60,0) {};
        \node[node] (B) at ( 2.60,0) {};
        \node[node] (C) at ( 0,1.90) {};
        \node[node] (D) at ( 0,0.43) {};

        \draw[thick] (A) -- (B);
        \draw[thick] (A) to[bend left=30] (C);
        \draw[thick] (A) -- (D);

        \draw[thick, ->, >=stealth] (A)
            to[bend right=14] node[label, circle, inner sep=0pt] {117} (B);
        \draw[thick, ->, >=stealth] (A)
            to[bend right=29] node[label, circle, inner sep=0pt] {147} (B);

        \draw[thick, ->, >=stealth] (B)
            to[bend right=30] node[label, circle, inner sep=0pt] {78} (C);
        \draw[thick, ->, >=stealth] (B)
            -- node[label, circle, inner sep=0pt] {12} (D);

        \draw[thick, ->, >=stealth] (C)
            to[bend right=65] node[label, pos=.45, circle, inner sep=0pt] {248} (D);
        \draw[thick, ->, >=stealth] (C)
            -- node[label, pos=.53, circle, inner sep=0pt] {173} (D);
        \draw[thick, ->, >=stealth] (C)
            to[bend left=65] node[label, pos=.60, circle, inner sep=0pt] {206} (D);
    \end{tikzpicture}
        \caption{$\texttt{SmallGroup}(264,12)=C_{11}\times\mathrm{SL}(2,3)$ and $n(5,9)\leq 1056$.}
        \label{fig:cageRecord59}
\end{subfigure}%

\caption{The base graphs and groups that resulted in a new upper bound for the cage problem. Darts corresponding to unmarked edges are assigned the neutral element. For Figs.~\ref{fig:cageRecord316},~\ref{fig:cageRecord317}~and~\ref{fig:cageRecord59}, group elements are numbered by their index as given by GAP's \texttt{MultiplicationTable}.}
\label{fig:liftBaseGraphs}
\end{figure}

Tables with all improved upper bounds are included in Appendix~\ref{app:tables}, while explicit constructions for the cage problem using lifts can be seen in Fig.~\ref{fig:liftBaseGraphs}. The previous arXiv version of this paper established $n(5,9) \leq 1116$. A subsequent search using similar voltage-lift methods yielded the construction in Fig.~\ref{fig:cageRecord59}, improving the bound to $n(5,9) \leq 1056$.

A selection of the currently best known bounds along with our improvements is available in Table~\ref{tab:cagecurrent1}, based on a survey of the literature. Lower bounds without a reference follow immediately from the Moore bound. We remark that in cases marked with a ``?'', it is known that a graph with the indicated property should exist, but the graph is orders of magnitude larger than the lower bound. Whenever the lower and upper bound coincide for a given parameter combination, the entire row is typeset in bold. The column labeled ``\# Graphs (this paper)'' indicates how many non-isomorphic graphs with those parameters were found using the algorithms described in this work. For the improved upper bounds obtained via excision, we halted the edge completion part (using the $(k,g)$-graph generator) after we found 100 graphs, while the generator could have generated more. That is why we indicate to have 100 pairwise non-isomorphic graphs for the excision cases, besides the improvement on $n(4,11)$ where we only obtained one graph.

An open conjecture for the cage problem states that every $(k,g)$-cage with even girth must be bipartite \cite{cagesurvey, wongcagesurvey}. The graphs associated with the new bounds for $n(3,16)$ and $n(4,10)$ are indeed bipartite, providing additional evidence for this conjecture.

For the studied variants, multiple new bounds were also obtained. In particular, we established 21 new bounds for $egr$-graphs and 29 for $vgr$-graphs, of which 2 and 7, respectively, match the corresponding lower bounds. The concrete bounds are provided in Tables~\ref{tab:edgeImproved} and \ref{tab:vertexImproved}. Additionally, for 34 $(k,g,n)$ tuples we determined the existence of a $(k,g)$-graph of order $n$ where it was previously not known whether $n$ is in the $(k,g)$-spectrum or not. These results are shown in Table~\ref{tab:spectraFound}. Finally, Table~\ref{tab:kgnoImproved} shows 6 improved upper bounds for $n(k,g,\uline{g+1})$. Applying the inequality $n(k,2s+4) \leq 2n(k,2s+1,\uline{2s+2})$ did not result in additional improvements for the cage problem.

To the best of our knowledge, upper bounds for $n(k,g,\lambda_v)$ with $k=4$ and $g \geq 7$, or with $k \in \{5,6\}$, were not available in the literature. While they are not listed in this work due to space constraints, they are available at \url{https://github.com/AGT-Kulak/smallRegGirthGraphs} for download.

\begin{table}[ht]
    \centering
    \caption{Best known lower and upper bounds for $n(k,g)$.}
    \label{tab:cagecurrent1}
\begin{tabular}{c c l l c c}
\toprule
&&& \multicolumn{2}{c}{$n(k,g) \leq$} &  \\ \cmidrule(lr){4-5}
        $k$ & $g$ & $n(k,g) \geq$ & Old & New & \# Graphs (this paper) \\\midrule
        \textbf{3} & \textbf{9} & \textbf{58 }\cite{Brinkmann1995TheSC} & \textbf{58} \cite{BIGGS1980299} & &\textbf{0} \\
        \textbf{3} & \textbf{10} & \textbf{70} \cite{WONG1983119} & \textbf{70} \cite{BALABAN19721, OKEEFE198091} & & \textbf{1}\\
        \textbf{3} & \textbf{11} & \textbf{112} \cite{10.5555/314613.314699}  & \textbf{112} \cite{balaban1973trivalent} & & \textbf{0}\\
        \textbf{3} & \textbf{12} & \textbf{126} & \textbf{126} \cite{Benson1966} & &\textbf{0}\\
        3 & 13 & 202 \cite{10.5555/314613.314699} & 272 \cite{biggs1989cubic} & & 1\\
        3 & 14 & 258 \cite{10.5555/314613.314699} & 384 \cite{exoosmalltriv} & & 1\\
        3 & 15 & 384 & 620 \cite{BiggsN.L.1983Tscf} & & 1\\
        3 & 16 & 512 & 960 \cite{exoovoltage2004} & 936& 1\\
        3 & 17 & 768 & 2176 \cite{exoo2019new} & 2048 & 1\\
        3 & 18 & 1024 & 2560 \cite{exoo2019new} & & 0\\
        3 & 19 & 1536 & 4324 \cite{HOAREM1993TaH} & & 0\\
        3 & 20 & 2048 & 5376 \cite{exoo2019new} & & 0\\
        \midrule

        \textbf{4} & \textbf{7} & \textbf{67} \cite{EXOO2011166} & \textbf{67} \cite{EXOO2011166} & & \textbf{0}\\
        \textbf{4} & \textbf{8} & \textbf{80} & \textbf{80} \cite{cagesurvey} & & \textbf{0}\\
        4 & 9 & 162 & 275 \cite{cagesurvey} & 270 & 2\\
        4 & 10 & 243 & 384 \cite{cagesurvey} & 320 & 1\\
        4 & 11 & 486 & ? & 713 & 1\\ 
        \textbf{4} & \textbf{12} & \textbf{728} & \textbf{728} \cite{cagesurvey} & & \textbf{0}\\
        \midrule

        5 & 7 & 108 & 152 \cite{cagesurvey} & & 0\\
        \textbf{5} & \textbf{8} & \textbf{170} & \textbf{170} \cite{cagesurvey} & & \textbf{0}\\
        5 & 9 & 428 & 1152 \cite{Exoo2025} & 1056 & 1\\
        5 & 10 & 684 & 1296 \cite{cagesurvey} & & 0\\
        5 & 11 & 1708 & 2688 \cite{ARAUJOPARDO20101622} & & 0\\
        \textbf{5} & \textbf{12} & \textbf{2730} & \textbf{2730} \cite{cagesurvey} & & \textbf{0}\\
        \midrule
        
        6 & 7 & 188 & 294 \cite{cagesurvey} & & 0\\
        \textbf{6} & \textbf{8} & \textbf{312} & \textbf{312} \cite{cagesurvey} & & \textbf{0}\\
        6 & 9 & 938 & ? & & 0\\
        6 & 10 & 1563 & ? & & 0\\
        6 & 11 & 4688 & ?  & 7783 & $\ge 100$ \\
        \textbf{6} & \textbf{12} & \textbf{7812} & \textbf{7812} \cite{cagesurvey} & & \textbf{0}\\
        \midrule

        8 & 7 & 457 & 777 \cite{deruiter2015applications} & 774 & $\ge 100$\\
        \midrule

        10 & 7 & 911 & 1611 \cite{deruiter2015applications} & 1608 & $\ge 100$\\
        \midrule

        12 & 7 & 1597 & 2893 \cite{deruiter2015applications} & 2890 & $\ge 100$\\
        \midrule

        14 & 7 & 2563 & 4719 \cite{deruiter2015applications} & 4716 & $\ge 100$\\
        
        \bottomrule
    \end{tabular}
\end{table}

\subsection{Validation and sanity checks}
In this subsection, we describe which measures we took to verify the correctness of our implementation of the algorithms described in Section~\ref{sec:algorithms}. Their correctness was verified empirically using multiple strategies.

First, the input was verified. For groups this is straightforward: all groups up to a given order are retrieved directly from GAP. For base graphs, the output of \texttt{multigraph+} was compared with the existing generator \texttt{pregraph} \cite{brinkmann2013generation} for cubic graphs containing loops and semi-edges. Both generators were tested up to order 13 and yielded identical results.

Second, potential implementation errors in the pruning optimizations using graph automorphisms $\phi_G$, group automorphisms $\phi_\Gamma$, and a girth bound $g_\text{min}$ were excluded by disabling them individually. In all cases exactly the same number of non-isomorphic lifts were recovered.

Finally, the exhaustiveness of the results was compared with independent implementations. BTA reproduced the same sets of non-isomorphic lifts as three external generators: the $K_{1,3}$-based construction from \cite{eze2024kggraphsg1cycles}, the Gray and Theta graph lifts from \cite{abreu2025graygraphpseudo2factor}, and a block-circulant graph generator \cite{GOEDGEBEUR2022212}, which is equivalent to lifts using cyclic groups. In each case the output was in complete agreement. 

Therefore, we are confident that our implementations are correct. Nevertheless, we stress that even if an implementation would be incorrect, the certificates that certify the bounds (i.e., the graphs that we found) are available for download (together with our source code) at \url{https://github.com/AGT-Kulak/smallRegGirthGraphs}. This makes it easier for others to independently verify our claims as it is not difficult to verify that the graphs have the properties that we claim.

\subsection{Algorithmic efficiency versus hardware scaling}

Given the long-standing nature of some of the previous bounds, it is natural to ask whether our success in finding new cages is merely the result of improvements in computing hardware. While access to modern high-performance computing infrastructure is obviously beneficial, our total computational effort was relatively modest and amounted to roughly five CPU-years. We therefore attribute our success primarily to the algorithmic setup. While the underlying theory of lifts of voltage graphs that we applied is not new, our implementation combining multiple pruning strategies and our approach of exhaustively and indiscriminately generating and testing all combinations of valid base graphs and groups up to a given order allowed us to locate candidate configurations that had been missed in earlier, more restricted searches.

To illustrate that the pruning techniques discussed in Section~\ref{sec:algorithms} are all of importance, Table~\ref{tab:sanity1} details the execution time of BTA under different combinations of exclusions. It is evident that the optimizations have vastly different impacts. For instance, pruning based on group automorphisms ($\phi_\Gamma$) or based on the minimum girth ($g_\text{min}$) yields massive performance gains individually, reducing the execution time by a factor of 100 or more. On the other hand, pruning based on edge automorphisms ($\phi_G$) provides a smaller speedup. More importantly, these optimizations combine together neatly, reducing the execution time by a factor of more than 30,000 when using all pruning rules compared to the unoptimized baseline. This shows that the new bounds are the result of algorithmic efficiency rather than just an increase in raw computing power.

\begin{table}[ht]
    \centering
    \caption{Execution time for exhaustively constructing lifts with BTA using 3-regular base graphs on three vertices, with groups up to order 40 and a minimum girth of 9, under various combinations of pruning optimizations. All versions found the same 191 graphs.}
    \label{tab:sanity1}
    \begin{tabular}{c c c S[table-format=5.1]} \toprule
        $\phi_G$ & $\phi_\Gamma$ & $g_\text{min}$ & {Execution time (s)} \\ \midrule
         & & & 16460.9\\
         \checkmark & & & 5849.9\\
         & \checkmark & & 160.0\\
         \checkmark& \checkmark & & 67.4\\\midrule
         &  & \checkmark& 123.4\\
         \checkmark&  & \checkmark& 19.7\\
         & \checkmark & \checkmark& 1.0\\
        \checkmark & \checkmark & \checkmark& 0.5\\
        \bottomrule
    \end{tabular}
\end{table}

\section{Conclusion}\label{sec:conclusion}
In this paper, we studied the cage problem and several of its variants from a computational perspective. Using three algorithms based on graph lifts and an excision strategy, we were able to establish multiple new upper bounds, including eleven for the classical cage problem. Notably, the upper bound for $n(4,10)$ was significantly improved from $384$ to $320$, while long-standing records for $n(3,16)$ and $n(3,17)$ were also broken. For the considered variants, our methods yielded dozens of new bounds, several of which match the corresponding lower bounds, and we determined for 34 tuples $(k,g,n)$ at least one $(k,g)$-graph of order $n$ for which it was previously not known whether $n$ is in the $(k,g)$-spectrum.

These results highlight the effectiveness of lift-based techniques for extremal graph theory problems. These newly discovered graphs may inspire new constructions leading to infinite families and we encourage further research in this direction. In addition, it would be interesting to systematically generalize the removed sets mentioned in the excision part for certain pairs $(k,g)$.

\section*{Acknowledgements}
 Jan Goedgebeur and Tibo Van den Eede are supported by Internal Funds of KU Leuven and a grant of the Research Foundation Flanders (FWO) with grant number G0AGX24N. Jorik Jooken is supported by an FWO Postdoctoral Fellowship with grant number 1222524N.

The computational resources and services used in this work were provided by the VSC (Flemish Supercomputer Center), funded by the Research Foundation - Flanders (FWO) and the Flemish Government.

\bibliographystyle{abbrvnatnourl}
\bibliography{references}
\clearpage

\appendix
\section{Tables}
\label{app:tables}

The following tables present analogous results for various variants of the cage problem considered in this work. The entry ``?'' indicates that no upper bound is known, bold rows indicate matching lower and upper bounds, and the ``\# Graphs'' column reports how many graphs were obtained by our algorithms.

Tables~\ref{tab:edgeImproved} and \ref{tab:vertexImproved} contain improved bounds for $egr$- and $vgr$-graphs respectively. All existing bounds for $vgr$-graphs follow from \cite{jajcay2024vgr}. Table~\ref{tab:spectraFound} contains previously undetermined $(k, g)$-spectra. Similarly, Table~\ref{tab:kgnoImproved} contains improved bounds for $(k, g)$-graphs without $(g+1)$-cycles.

All graphs found through lifts or excision are available at \url{https://github.com/AGT-Kulak/smallRegGirthGraphs}. The graphs for the cage problem are also available through \href{https://houseofgraphs.org/meta-directory/cages}{The House of Graphs} \cite{houseofgraphs} at \url{https://houseofgraphs.org/meta-directory/cages}.

\begin{table}[ht]
\centering
\caption{Overview of improved $n(k, g, \lambda_e)$-bounds.}
\label{tab:edgeImproved}
\centering
\begin{tabular}{c c c c c c c}
\toprule
&&&&\multicolumn{2}{c}{$n(k,g,\lambda_e) \leq$} & \\ \cmidrule(lr){5-6}
$k$ & $g$ & $\lambda_e$ & $n(k,g,\lambda_e) \geq$ & Old & New & \# Graphs (this paper) \\ \midrule
\textbf{4} & \textbf{5} & \textbf{3} & \textbf{30} \cite{goedgebeur2024egrgraphs} & \textbf{55} \cite{goedgebeur2024egrgraphs} & \textbf{30} & \textbf{1} \\
\textbf{4} & \textbf{5} & \textbf{4} & \textbf{30} \cite{goedgebeur2024egrgraphs} & \textbf{?} & \textbf{30} & \textbf{4}  \\
\midrule

4 & 6 & 1 & 57 \cite{goedgebeur2024egrgraphs} & 84 \cite{goedgebeur2024egrgraphs} & 72 & 2 \\
4 & 6 & 2 & 51 \cite{goedgebeur2024egrgraphs} & 96 \cite{goedgebeur2024egrgraphs} & 69 & 1 \\
4 & 6 & 3 & 45 \cite{goedgebeur2024egrgraphs} & 60 \cite{goedgebeur2024egrgraphs} & 57 & 1 \\
4 & 6 & 5 & 42 \cite{goedgebeur2024egrgraphs} & 81 \cite{goedgebeur2024egrgraphs} & 60 & 3 \\
4 & 6 & 6 & 40 \cite{goedgebeur2024egrgraphs} & 64 \cite{goedgebeur2024egrgraphs} & 48 & 2 \\
4 & 6 & 7 & 39 \cite{goedgebeur2024egrgraphs} & 60 \cite{goedgebeur2024egrgraphs} & 54 & 1 \\
4 & 6 & 15 & 33 \cite{goedgebeur2024egrgraphs} & 40 \cite{goedgebeur2024egrgraphs} & 35 & 1 \\
\midrule

5 & 5 & 2 & 44 \cite{goedgebeur2024egrgraphs} & 60 \cite{goedgebeur2024egrgraphs} & 56 & 1 \\
5 & 5 & 4 & 38 \cite{DrglinAjdaZavrtanik2021EEG} & 66 \cite{goedgebeur2024egrgraphs} & 48 & 1 \\
5 & 5 & 6 & 36 \cite{DrglinAjdaZavrtanik2021EEG} & ? & 44 & 1 \\
\midrule

6 & 4 & 7 & 24 \cite{goedgebeur2024egrgraphs} & 36 \cite{goedgebeur2024egrgraphs} & 28 & 1 \\
\midrule

6 & 5 & 1 & 65 \cite{DrglinAjdaZavrtanik2021EEG} & ? & 120 & 1 \\
6 & 5 & 2 & 60 \cite{DrglinAjdaZavrtanik2021EEG} & 7620 \cite{JAJCAY201870} & 160 & 1 \\
6 & 5 & 4 & 60 \cite{DrglinAjdaZavrtanik2021EEG} & 910 \cite{JAJCAY201870} & 110 & 9 \\
6 & 5 & 5 & 57 \cite{DrglinAjdaZavrtanik2021EEG} & ? & 80 & 1 \\
6 & 5 & 6 & 60 \cite{DrglinAjdaZavrtanik2021EEG} & ? & 100 & 2 \\
6 & 5 & 7 & 55 \cite{DrglinAjdaZavrtanik2021EEG} & ? & 100 & 1 \\
6 & 5 & 9 & 55 \cite{DrglinAjdaZavrtanik2021EEG} & ? & 155 & 1 \\
6 & 5 & 10 & 52 \cite{DrglinAjdaZavrtanik2021EEG} & ? & 64 & 9 \\
\bottomrule
\end{tabular}
\end{table}

\begin{table}
\centering
\caption{Overview of improved $n(k, g, \lambda_v)$-bounds.}
\label{tab:vertexImproved}
\centering
\begin{tabular}{c c c c c c c}
\toprule
&&&&\multicolumn{2}{c}{$n(k,g,\lambda_v) \leq$} &  \\ \cmidrule(lr){5-6}
$k$ & $g$ & $\lambda_v$ & $n(k,g,\lambda_v) \geq$ & Old & New & \# Graphs (this paper) \\ \midrule
\textbf{3} & \textbf{7} & \textbf{1} & \textbf{42} & \textbf{56} & \textbf{42} & \textbf{2} \\
\textbf{3} & \textbf{7} & \textbf{2} & \textbf{42} & \textbf{?} & \textbf{42} & \textbf{3} \\
\textbf{3} & \textbf{7} & \textbf{4} & \textbf{42} & \textbf{?} & \textbf{42} & \textbf{2} \\
\midrule

3 & 8 & 2 & 52 & 64 & 60 & 3 \\
\textbf{3} & \textbf{8} & \textbf{5} & \textbf{48} & \textbf{64} & \textbf{48} & \textbf{1} \\
3 & 8 & 7 & 48 & ? & 56 & 3 \\
\midrule

\textbf{4} & \textbf{5} & \textbf{1} & \textbf{35} & \textbf{420} & \textbf{35} & \textbf{3}  \\
4 & 5 & 3 & 30 & 40 & 35 & 5  \\
\textbf{4} & \textbf{5} & \textbf{6} & \textbf{30} & \textbf{55} & \textbf{30} & \textbf{22} \\
\textbf{4} & \textbf{5} & \textbf{7} & \textbf{30} & \textbf{?} & \textbf{30} & \textbf{9} \\
\midrule

4 & 6 & 1 & 54 & 1152 & 72 & 3 \\
4 & 6 & 2 & 51 & 84 & 72 & 8 \\
4 & 6 & 3 & 50 & 1152 & 60 & 1 \\
4 & 6 & 4 & 48 & 96 & 60 & 1 \\
4 & 6 & 5 & 48 & ? & 60 & 2 \\
4 & 6 & 6 & 46 & 60 & 55 & 1 \\
4 & 6 & 7 & 48 & 648 & 60 & 16 \\
4 & 6 & 8 & 45 & 60 & 54 & 2 \\
4 & 6 & 9 & 44 & 512 & 48 & 1 \\
4 & 6 & 10 & 42 & 81 & 51 & 2 \\
4 & 6 & 11 & 42 & ? & 54 & 3 \\
4 & 6 & 12 & 40 & 64 & 48 & 3 \\
4 & 6 & 13 & 42 & ? & 48 & 2 \\
4 & 6 & 14 & 39 & 60 & 54 & 11 \\
4 & 6 & 15 & 38 & ? & 48 & 1 \\
4 & 6 & 17 & 42 & ? & 54 & 18 \\
4 & 6 & 19 & 36 & ? & 48 & 9 \\
4 & 6 & 21 & 36 & 44 & 42 & 2 \\
4 & 6 & 29 & 36 & ? & 42 & 2 \\
\bottomrule
\end{tabular}
\end{table}

\begin{table}[ht]
\centering
\caption{Overview of found $(k, g)$-graphs of order $n$ for which it was open whether $n$ is in the $(k,g)$-spectrum as described in \cite{eze2025spectra}.}
\label{tab:spectraFound}

\makebox[\linewidth][c]{
\hspace{4em}
\begin{subtable}{0.25\linewidth}
\centering
\caption{The cases $k=3$ and $k=4$.}
\begin{tabular}{c c c c}
\toprule
$k$ & $g$ & $n$ & \# Graphs (this paper) \\ \midrule
3 & 11 & 140 & 3  \\
\midrule
3 & 12 & 152 & 1  \\
3 & 12 & 176 & 1  \\
3 & 12 & 186 & 1  \\
3 & 12 & 190 & 1  \\
3 & 12 & 198 & 9  \\
3 & 12 & 200 & 15 \\
\midrule
4 & 7 & 81 & 3 \\
\midrule
4 & 8 & 90  & 1 \\
4 & 8 & 112 & 3 \\
4 & 8 & 117 & 4 \\
4 & 8 & 123 & 1 \\
4 & 8 & 129 & 3 \\
\bottomrule
\end{tabular}
\end{subtable}
\hspace{1em}
\bigskip
\begin{subtable}{0.65\linewidth}
\centering
\caption{The case $k=5$.}
\begin{tabular}{c c c c}
\toprule
$k$ & $g$ & $n$ & \# Graphs (this paper) \\ \midrule
5 & 7 & 184 & 2 \\
5 & 7 & 188 & 7 \\
5 & 7 & 190 & 32 \\
5 & 7 & 192 & 69 \\
5 & 7 & 198 & 5 \\
5 & 7 & 200 & 7 \\
5 & 7 & 204 & 2 \\
5 & 7 & 208 & 42 \\
5 & 7 & 210 & 18 \\
5 & 7 & 212 & 19 \\
5 & 7 & 216 & 253 \\
5 & 7 & 220 & 258 \\
5 & 7 & 224 & 943 \\
5 & 7 & 228 & 152 \\
5 & 7 & 230 & 5185 \\
5 & 7 & 232 & 217 \\
5 & 7 & 234 & 41 \\
5 & 7 & 236 & 453 \\
5 & 7 & 240 & 10703 \\
5 & 7 & 244 & 745 \\
5 & 7 & 248 & 891 \\
5 & 7 & 250 & 46 \\
5 & 7 & 252 & 214 \\
5 & 7 & 256 & 2531 \\
\bottomrule
\end{tabular}
\end{subtable}
}
\end{table}

\begin{table}[ht]
\caption{Overview of improved bounds for ($k, g$)-graphs without ($g+1$)-cycles for odd $g$.}
\label{tab:kgnoImproved}
\centering
\begin{tabular}{c c c l c c}
\toprule
&&& \multicolumn{2}{c}{$n(k,g,\uline{g+1}) \leq$} & \\ \cmidrule(lr){4-5}
$k$ & $g$ & $n(k,g,\uline{g+1}) \geq$ & Old & New  & \# Graphs (this paper) \\ \midrule
3 & 11 & 144 \cite{eze2024kggraphsg1cycles} & 272 \cite{potocnik2012cubicvertextransitivegraphs1280} & 228 & 1 \\
3 & 13 & 286 \cite{eze2024kggraphsg1cycles} & 800 \cite{potocnik2012cubicvertextransitivegraphs1280} & 600 & 1\\
3 & 15 & 576 \cite{eze2024kggraphsg1cycles} & 2162 \cite{condermarstontriv} & 1760 & 3 \\
\midrule

4 & 7 & 126 \cite{eze2024kggraphsg1cycles} & 273 \cite{potocnik2012cubicvertextransitivegraphs1280, POTOCNIK2015148} & 195  & 1 \\
4 & 9 & 377 \cite{eze2024kggraphsg1cycles} & 1518 \cite{POTOCNIK20091323} & 972 & 1 \\

\midrule
5 & 5 & 86 \cite{eze2024kggraphsg1cycles} & 192 \cite{potocnik2024census} & 128 & 3 \\
\bottomrule
\end{tabular}
\end{table}

\end{document}